\documentclass[10pt]{article} 
\usepackage[utf8]{inputenc}
\usepackage[T1]{fontenc}
\usepackage{lmodern}
\usepackage{amsmath}
\usepackage{amssymb}
\usepackage{amsthm}
\usepackage{graphics}
\usepackage{graphicx}
\usepackage{stmaryrd}
\usepackage{xcolor,colortbl}
\definecolor{mycolor}{rgb}{.8, .8, 1}
\usepackage{caption}
\usepackage{placeins}
\usepackage{bbm}
\usepackage{multirow}
\usepackage{booktabs}
\usepackage{comment}
\usepackage{utfsym}
\usepackage{tikz}

\usepackage{subcaption}
\usepackage[ruled]{algorithm2e}
\usepackage{setspace} 

\usepackage{todonotes}

\usepackage{hyperref}
\hypersetup{
    colorlinks,
    citecolor=black,
    filecolor=black,
    linkcolor=black,
    urlcolor=black
}

\usepackage[top=1.5in, bottom=1.5in, left=1.5in, right=1.5in]{geometry}

\usepackage{xurl} 
\usepackage[capitalize,nameinlink]{cleveref}

\crefname{figure}{fig.}{figs.}
\Crefname{figure}{Fig.}{Figs.}

\crefname{table}{tab.}{tabs.}
\Crefname{table}{Tab.}{Tabs.}

\crefname{equation}{eq.}{eqs.}
\Crefname{equation}{Eq.}{Eqs.}

\crefname{section}{sec.}{secs.}
\Crefname{section}{Sec.}{Secs.}

\crefname{appendix}{app.}{apps.}
\Crefname{appendix}{Appendix}{Appendix}

\crefname{property}{prop.}{props.}
\Crefname{property}{Prop.}{Props.}

\crefname{definition}{def.}{defs.}
\Crefname{definition}{Def.}{Defs.}

\crefname{remark}{rem.}{rems.}
\Crefname{remark}{Rem.}{Rems.}

\crefname{theorem}{thm.}{thms.}
\Crefname{theorem}{Thm.}{Thms.}

\newcommand{\sci}[1]{\!\times\! 10^{#1}}

\def\R{ { \mathbb{R} } }

\def\E{ { \mathbf{E} } }

\def \X{ { x } } 
\def \V{ { v } } 
\def \U{ { u } } 

\def \bfX{ { \mathbf{\X}  } }
\def \bfV{ { \mathbf{\V}  } }
\def \bfU{ { \mathbf{\U}  } }

\def \x{ { \bar{x} } }

\def \bfx{ { \mathbf{\x}  } }
\def \bfv{ { \mathbf{\bar{v} }  } }
\def \bfu{ { \mathbf{\bar{u}}  } }

\def \hami{ { \mathcal{H}  } }
\def \calD{ { \mathcal{D}  } }
\def \calE{ { \mathcal{E}  } }

\def \half{ { \tfrac{1}{2} } }

\def \opmean{\operatorname{mean}}
\def \opmin{\operatorname{min}}
\def \opmax{\operatorname{max}}

\def \AE{\text{AE}}

\usepackage [mathscr] {eucal}
\def \ee{\frac{1}{2}\|E\|_{2}}
\def \error{\text{err}}

\def \model{\text{pred}}

\def \modelr{{\text{p}\overline{\text{re}}\text{d}}}

\def \stab{{\text{s}\overline{\text{ta}}\text{b}}}

\def \init{\text{init}}
\def \train{\text{train}}
\def \test{\text{test}}
\def \reference{\text{ref}}

\newtheorem*{remark}{Remark}

\title{Reduced Particle in Cell method for the Vlasov-Poisson system using auto-encoder and Hamiltonian neural networks}

\usepackage{authblk}
\author[1]{Raphaël Côte}
\author[1, 2]{Emmanuel Franck}
\author[1, 2]{Laurent Navoret}
\author[1, 2]{Guillaume Steimer}
\author[1, 2]{Vincent Vigon}
\affil[1]{Institut de Recherche Mathématique Avancée, UMR 7501 Université de Strasbourg et CNRS, 7 rue René Descartes 67000 Strasbourg, France}
\affil[2]{INRIA Nancy-Grand Est, MACARON Project, Strasbourg, France}

\date{}

\usepackage[backend=bibtex]{biblatex}
\begin{document}

\sloppy
\maketitle 

\begin{abstract}
    Hamiltonian particle-based simulations of plasma dynamics are inherently computationally intensive, primarily due to the large number of particles required to obtain accurate solutions. This challenge becomes even more acute in many-query contexts, where numerous simulations must be conducted across a range of time and parameter values. Consequently, it is essential to construct reduced order models from such discretizations to significantly lower computational costs while ensuring validity across the specified time and parameter domains. Preserving the Hamiltonian structure in these reduced models is also crucial, as it helps maintain long-term stability. 
    
    In this paper, we introduce a nonlinear, non-intrusive, data-driven model order reduction method for the 1D-1V Vlasov–Poisson system, discretized using a Hamiltonian Particle-In-Cell  scheme. Our approach relies on a two-step projection framework: an initial linear projection based on the Proper Symplectic Decomposition, followed by a nonlinear projection learned via an autoencoder neural network. The reduced dynamics are then modeled using a Hamiltonian neural network. The offline phase of the method is split into two stages: first, constructing the linear projection using full-order model snapshots; second, jointly training the autoencoder and the Hamiltonian neural network  to simultaneously learn the encoder-decoder mappings and the reduced dynamics. We validate the proposed method on several benchmarks, including Landau damping and two-stream instability. The results show that our method has better reduction properties than standard linear Hamiltonian reduction methods.
\end{abstract}
\vspace{2em}

\\


\clearpage

\section{Introduction}

Plasma are gases made of charged particles interacting through long-range Coulomb interactions.  A standard kinetic approach for characterizing collisionless plasma dynamics is based on the Vlasov–Maxwell equations, which describe the time evolution of  the particle distribution functions in position-velocity phase space and the dynamics of the self-consistent electromagnetic fields. In this work, we focus on the electrostatic limit of these equations, namely the Vlasov–Poisson system, where particle interactions are driven by a self-consistent electric field satisfying the Poisson equation.

Simulating the Vlasov–Poisson system numerically presents significant challenges, and a wide range of particle-based methods have been developed to address them. These methods represent the charged particles distribution using a large set of macro particles, whose trajectories are evolved according to the characteristics of the kinetic equation. To compute the self-induced electric field, the computational domain is discretized into a mesh on which the Poisson equation is solved. The particle distribution is projected onto this mesh to get the charge density, the electric field is computed there, and then interpolated back to the particle positions. The particles are subsequently moved under the influence of the resulting Lorentz force. This approach is known as the Particle-In-Cell (PIC) method \cite{birdsall2004plasma, pritchett2003particle}.

Over time, PIC methods have been refined to preserve important physical invariants, such as total energy (see \cite{lewis1970energy,chen2011energy} and references therein). Notably, the Vlasov–Poisson system admits a Hamiltonian formulation \cite{marsden1982hamiltonian}, which ensures conservation of total energy and the system’s underlying symplectic structure. Preserving this Hamiltonian structure in the discretized PIC framework is essential for maintaining long-term numerical stability. As a consequence, several Hamiltonian PIC methods have been developed, including the Hamiltonian PIC scheme \cite{he2016hamiltonian}, as well as canonical and non-canonical symplectic PIC methods \cite{qin2016canonical, xiao2015explicit}, and the Geometric Electromagnetic PIC (GEMPIC) method \cite{kraus2017gempic}.

\vspace{1em}

Given the nonlinear dynamics and multiscale phenomena of the Vlasov–Poisson system, along with the need to employ a very large number of particles to achieve accurate convergence to the solution \cite{barsamian2018picvert}, PIC simulations present a significant numerical challenge. This makes the use of Hamiltonian model order reduction techniques particularly compelling. In real time or many query contexts—such as control processes, optimization, or uncertainty quantification—reduced order models (ROMs) can be crucial. Starting from a particle-based discretization of the Vlasov–Poisson system, referred to as the full order model (FOM), model reduction seeks to construct a smaller dynamical system that provides accurate approximations over a specified range of times and parameters, while substantially lowering computational cost. Crucially, preserving the Hamiltonian structure in the reduced model contributes to its long-term robustness and stability \cite{tyranowski2023symplectic}.

Over the recent years, considerable efforts have been devoted to constructing reduced models for the Vlasov–Poisson dynamics. These surrogate models aim to reduce the computational cost associated with plasma simulations. By allowing for a small approximation error, reduced models enable significantly faster solution evaluations. Broadly speaking, three main families of model order reduction techniques have been applied to the Vlasov–Poisson system. These approaches are well reviewed in \cite{nouy2017lowrank} within a more general framework.

The first family comprises projection-based model order reduction methods \cite{grassle2021model, hesthaven2022reduced, hesthaven2024adaptive}. These approaches assume that the solution manifold lies close to a low-dimensional subspace, which is approximated using a collection of solution snapshots obtained thanks to Proper Orthogonal Decomposition (POD) or greedy approach. The reduced model is then derived using a Galerkin projection onto this subspace. However, projection alone often does not suffice to reduce computational cost: to compute the reduced system’s vector field, the reduced state must typically be lifted back to the full physical space. To mitigate this costly round trip, various strategies have been introduced. For example, the Discrete Empirical Interpolation Method (DEIM) \cite{chaturantabut2010nonlinear} selects a small set of spatial points at which to evaluate the nonlinear terms, while Dynamic Mode Decomposition (DMD) \cite{tissot2014model, schmid2010dynamic} extracts relevant modes from the data and models their evolution with a simplified linear system.
The second family concerns sparse approximation methods, in which solutions are approximated by selecting a small number of basis functions, based on prior knowledge about the structure of the solution. This has been explored for the Vlasov-Poisson dynamics in \cite{kormann2016sparse}, where the dynamics is computed on a sparse grid using a semi-Lagrangian solver. 
The third family is made of low-rank approximation methods, where the solutions are expressed as a sum of low rank tensors. This has been applied for plasma dynamics in \cite{ehrlacher2017dynamical, einkemmer2021mass, einkemmer2018lowrank}. For instance, \cite{einkemmer2021mass} proposes a continuous low-rank representation of the distribution function, followed by discretization using a conservative dynamical low-rank scheme that preserves key physical quantities such as mass, momentum, and energy.
\vspace{1em}

While these techniques have proven effective for continuous problems in Eulerian or semi-Lagrangian frameworks, the model reduction of particle-based discretizations of the Vlasov–Poisson system poses a distinct challenge. In parallel with the development of these methods, several studies have highlighted the potential of machine learning—particularly neural networks—to assist with or even automate aspects of model order reduction. Examples include manifold learning techniques \cite{hyperelastic2016nonlinear,sonday2010manifold}, as well as data-driven reduced order models for PDEs derived from mesh-based discretizations \cite{kim2022fast, lee2020model}.

In addition, preserving the Hamiltonian structure in the reduced model is a supplementary challenge. Using reduction techniques which do not preserve structure, such as POD, often leads to numerically unstable models: their dynamics can diverge significantly from the underlying physical behavior. Fortunately, several reduction techniques can be adapted to retain structural properties within the reduced model. For example, the DEIM method can be modified to preserve the first moments of an operator, as shown in \cite{franck2024hyperbolic}. Specifically for Hamiltonian systems, structure preservation can be ensured through the Proper Symplectic Decomposition (PSD) \cite{peng2015symplectic}, a symplectic counterpart to POD, in which the projection onto the reduced space is constrained to be symplectic. For the PIC model that motivates our work, \cite{hesthaven2024adaptive} introduces a dynamic, projection-based model order reduction framework. The projection evolves in time and, with additional constraints, can be made symplectic—thus ensuring that the reduced model remains Hamiltonian. To further enhance computational efficiency, their approach also integrates a DMD-DEIM method to reduce the number of particles effectively. Machine learning techniques have also been extended to account for Hamiltonian structures. In \cite{buchfink2023symplectic}, for instance, the authors replace the linear PSD mapping with a neural network that is weakly constrained to be symplectic. 

\vspace{1em}
 
In this paper, we focus on the 1D-1V Vlasov–Poisson system discretized using a Hamiltonian PIC scheme. This discretization yields a high-dimensional ODE with a Hamiltonian structure. However, relying solely on a PSD to construct a reduced model is insufficient to capture small-scale dynamics and nonlinear behaviors with a small reduced dimension. As a result, such an approach would offer limited computational speedup.

To address this issue, we present a strategy inspired by \cite{fresca2022poddlrom}, where the authors perform an initial projection onto an intermediate subspace using a Proper Orthogonal Decomposition (POD), followed by the construction of a reduced model through deep learning techniques. Our approach, which we refer to as the PSD-AE-HNN method, similarly employs a two-step projection. It combines PSD with the AE-HNN method introduced in \cite{cote2025hamiltonian}.

Starting from the full state variables of the PIC model, we first apply a PSD-based projection onto a symplectic subspace  of intermediate dimension. Due to the symplectic nature of the PSD, the intermediate variables follows a Hamiltonian dynamic. Next, we perform a second, nonlinear projection using an autoencoder (AE) neural network \cite{goodfellow2016deep} to map the system onto a lower dimensional space. While this second mapping is not explicitly symplectic, we enforce Hamiltonian structure in the reduced dynamics by training a Hamiltonian Neural Network (HNN) \cite{greydanus2019hamiltonian} and incorporating tailored loss functions to constrain the training process.

The motivation behind this two-step mapping lies in the nature of the PIC discretization: particles are neither ordered nor regularly spaced, precluding the use of convolutional neural networks. Furthermore, a large number of particles (e.g. $10^5$ in 1D–1V) are typically required to achieve accurate convergence, making the direct use of dense neural networks impractical. The PSD thus acts as a symplectic preconditioner, enabling the AE-HNN method to learn dynamics from a Hamiltonian intermediate representation of reasonable size (e.g. $10^2$).

The structure of this paper is as follows: In the first section, we recall the Vlasov-Poisson equation and its Hamiltonian PIC discretization. In the second section, we present our model order reduction technique, referred to as the PSD-AE-HNN method. In the third section, we apply this method to several classic numerical test cases, including linear and nonlinear Landau damping and the two-stream instability. We then provide a comparison of computational times before concluding.

\section{Particle discretization of the Vlasov-Poisson equation}

In this section, we present the full order model (FOM) of interest. It is a particle-based discretization of the Vlasov-Poisson equation which possesses a Hamiltonian structure. 

    \subsection{The Vlasov-Poisson equation}

We consider a parametric 1D-1V Vlasov-Poisson equation, that gives the dynamics of the particle distribution function $f(t, x, v; \mu)$ which depends on time $t \in [0, T]$ with $T>0$, position $x \in \Omega_x$, a periodic domain of size $|\Omega_x |$, velocity $v \in \Omega_v \subset \R$, and parameters $\mu \in \Gamma \subset \R^d$ with $d > 0$. The equation reads
\begin{align}\label{eq:vlasov-poisson}
    \begin{cases}
        \displaystyle\partial_t f(t, x, v; \mu) + v \, \partial_x f(t, x, v; \mu) + \frac{q}{m} E (t,x;\mu) \, \partial_v f(t, x, v; \mu) = 0, &\text{ in } [0,T] \times  \Omega_x\times\Omega_v\times\Gamma, \\
         \displaystyle\partial_{x} E (t,x;\mu) = q \int_{\Omega_v} f(t, x, v; \mu) \, dv - \rho_0, &\text{ in } [0,T] \times \Omega_x\times\Gamma, \\
        f(0,x,v;\mu) = f_{\init}(x,v;\mu), &\text{ in } \Omega_x\times\Omega_v\times\Gamma,
    \end{cases}
\end{align}
where $E(t, x; \mu) \in \R$ is the electric field, $q$ is the individual charge of the particles, $m$ their individual mass and $f_{\init}(x,v;\mu) \in \R$ is a given initial condition. Defining the charge density $\rho(t, x; \mu) = q \int_\R f(t, x, v; \mu) \, dv$ and the electric potential $\phi (t, x; \mu) \in \R$ such that $E(t, x; \mu) = - \partial_x \phi(t, x; \mu)$, the Poisson equation rewrites
\begin{equation}\label{eq:poisson from vlasov}
    - \partial_{xx} \phi (t,x;\mu) = \rho(t, x; \mu) - \rho_0, \quad \text{ in } [0,T] \times \Omega_x\times\Gamma.
\end{equation}
The variable $\rho_0$ corresponds to the average global charge density initially and remains constant over time:
\begin{equation}
     \rho_0 = \frac{1}{|\Omega_x|}\int_{\Omega_x} \rho(t,x;\mu) \, dx .\label{eq:normalization}
\end{equation} This quantity is subtracted from the charge density $\rho$ in the right-hand side of the Poisson equation to ensure that the system is well posed with periodic boundary conditions. 

The Vlasov-Poisson equation given in \Cref{eq:vlasov-poisson} admits a Hamiltonian formulation with a Lie-Poisson bracket \cite{casas2017high}, and a Hamiltonian which corresponds to the sum of the kinetic and potential energies of the system 
\begin{equation*}
    H(f;\mu) = \frac{m}{2}
    \int_{\Omega_x\times\Omega_v} v^2 f(t,x,v;\mu) \, dx dv
    + \frac{1}{2} \int_{\Omega_x}
    \left| E(t,x;\mu) \right|^2 \, dx.
\end{equation*}

    \subsection{Hamiltonian particle-based discretization}\label{sec : PIC disretization}

We consider a Particle-In-Cell (PIC) discretization of \Cref{eq:vlasov-poisson} that preserves the Hamiltonian structure of the equations. Namely, we use a particularization of the GEMPIC algorithm \cite{kraus2017gempic}, as considered in \cite{hesthaven2024adaptive}. The distribution function $f$ is approximated with a set of $N \in \mathbb{N}$ macro-particles and the electric field is obtained by solving the Poisson equation with a finite element discretization resulting in an approximate electric field $E_h(t,x,\mu)$. More precisely, we approximate $f$ with a sum of Dirac delta distributions, located at  position $(x_k(t;\mu) ,v_k(t;\mu))$ in phase space: 
\begin{equation*}
    f_N (t, x, v; \mu) = \sum_{k=1}^N \omega \, \delta \left(x - x_k(t;\mu)\right) \delta \left(v - v_k(t;\mu)\right), 
\end{equation*}
where $\omega$ is the weight of each particle assumed to be identical for all particles and set equals to $|\Omega_x|\rho_0/(qN)$ to ensure the charge density to be normalized (\Cref{eq:normalization}). To satisfy the Vlasov equation of \Cref{eq:vlasov-poisson}, the dynamics of $N$ particles have to satisfy the following system of differential equations:
\begin{align}
    \begin{cases}\label{eq:vlasov-poisson particle}
        \displaystyle\dfrac{d}{dt} \bfX(t;\mu) = \bfV(t;\mu), &\text{ in } [0,T], \\[10pt]
        \displaystyle\dfrac{d}{dt} \bfV(t;\mu) = \frac{q}{m} \E_h(t, \bfX(t;\mu);\mu), &\text{ in } [0,T], \\[10pt]
        \bfX(0;\mu) = \bfX_{\init}(\mu), \\
        \bfV(0;\mu) = \bfV_{\init}(\mu),
    \end{cases}
\end{align}
where $\bfX(t;\mu) = (x_k(t;\mu)), \bfV(t;\mu) = (v_k(t;\mu))\in \R^{N}$ denotes the vectors of positions and velocities and $\E_h(t, \bfX;\mu) = (E_h(t, x_k; \mu)) \in \R^N$ the approximate electric field evaluated at each particle position. 
To obtain this approximate electric field, an approximate charge density $\rho_h(t,x;\mu)$ is computed on the finite element mesh from the particles distribution (deposition step), the Poisson equation is solved and then the electric field is evaluated at particle positions (interpolation step). We thus need deposition and interpolation steps such that the resulting system is still Hamiltonian.  

In detail, we introduce a uniform grid of $\Omega_x$, denoted $X_h=\left\{ih, i \in  \left\{1, \cdots, n_x\right\} \right\}$, where $h$ is the cell length. We consider a $H^1$-conforming finite element discretization of the Poisson \Cref{eq:poisson from vlasov} in the space $\mathcal{P}_1 \Lambda^0 (\Omega_x)$ of piecewise linear functions. As in \cite{hesthaven2024adaptive}, $\left( \lambda^0_i (x) \right)_{i \in \{1,\cdots,n_x\}}$ denotes the basis, which satisfies $\lambda^0_i(jh)=\delta_{i,j}$ with $\delta_{i,j}$ the Kronecker delta. Then, we define the particle-to-grid mapping $\Lambda^0(\bfX) \in \mathcal{M}_{N,n_x}(\R)$:
\[
    \left( \Lambda^0(\bfX) \right)_{k,i} = \lambda^0_i(x_k), \qquad k \in \{1, \cdots, N\}, i \in \{1, \cdots, n_x\},
\]
and the matrix of the Poisson problem $L \in \mathcal{M}_{n_x,n_x}(\R)$ by
\[
    L_{i,j} = \langle d_x \lambda^0_i, d_x \lambda^0_j \rangle_{L^2(\Omega_x)}, \qquad i,j \in \{1, \cdots, n_x\},
\]
where $d_x$ is the derivative with respect to $x$ and $\langle \cdot, \cdot \rangle_{L^2(\Omega_x)}$ is the standard $L^2(\Omega_x)$ scalar product, which in our case equals the standard one-dimensional discrete Laplacian matrix (up to factor $1/h$):
\begin{equation*}
    L= \frac{1}{h} 
        \begin{pmatrix}
            -2 & 1 &  &  \\
            1 & \ddots & \ddots &  \\
             & \ddots & \ddots & 1 \\
             &  & 1 & -2 
        \end{pmatrix}.
\end{equation*}
With these notations, the computation of the approximate electric field can be written as follows. From the particles positions, we compute a discrete charge density
    \begin{align*}
        &\text{ (deposition step) }&&\boldsymbol{\rho}_h = q\omega \, \Lambda^0(\bfX)^T \mathbbm{1}_N = \left(q\omega \sum_{k=1}^N \lambda^0_i(x_k)\right)_i \in \R^{n_x}.\\
        \intertext{where $\mathbbm{1}_N \in \R^N$ is a vector of ones. Then, the approximate potential is computed by solving the discrete Poisson equation}
        &&& - L \boldsymbol{\phi}_h = \boldsymbol{\rho}_h - h \rho_0\mathbbm{1}_{n_x},\\
        \intertext{with $\boldsymbol{\phi}_h \in \R^{n_x}$. Finally, the discrete electric field is defined by $E_h(t,x;\mu) = -\sum_{i=1}^{n_x} d_x \lambda^0_i(x) (\boldsymbol{\phi}_h)_i$, and can be evaluated at particles positions:}
        &\text{ (interpolation step)}&& \mathbf{E}_h = -\nabla \Lambda^0\left(\bfX\right)\boldsymbol{\phi}_h = \left(-\sum_{i=1}^{n_x} d_x \lambda^0_i(x_k) (\boldsymbol{\phi}_h)_i\right)_k \in \R^N,
    \end{align*}
with $\nabla \Lambda^0(\bfX) = (d_x\lambda^0_i(x_k))_{k,i}\in\mathcal{M}_{N,n_x}(\R)$. We note that we recover the deposition and interpolation steps of the standard PIC method \cite{birdsall2004plasma}.

The resulting system has a Hamiltonian structure. Indeed, introducing the discrete Hamiltonian function
\begin{align}\label{eq:H vlasov poisson}
    &\hami \left(\bfX(t;\mu),\bfV(t;\mu)\right) = \frac{1}{2} \|\bfV(t;\mu)\|^2 + \mathcal{U}\left( \bfX(t;\mu)\right), \\
    \intertext{with }
    &\mathcal{U}\left( \bfX(t;\mu) \right)
    = \frac{1}{2m\omega}     \left(q\omega \Lambda^0\left(\bfX(t;\mu)\right)^T\mathbbm{1}_N - h \rho_0\mathbbm{1}_{n_x}\right)^T L^{-1} \left( q\omega\Lambda^0\left(\bfX(t;\mu)\right)^T \mathbbm{1}_N - h \rho_0\mathbbm{1}_{n_x} \right),  
\end{align}
and the variable $\bfU(t;\mu) = (\bfX(t;\mu),\bfV(t;\mu))$, the full order dynamics \Cref{eq:vlasov-poisson particle} rewrite as a Hamiltonian system:
\begin{align}\label{eq:hamiltonian fom}
    \begin{cases}
    \displaystyle\frac{d}{dt} \bfU(t;\mu) = J_{2N} \nabla_\bfU \hami\left(\bfU(t;\mu)\right) , &\text{ in } [0,T]\\
    \bfU(0;\mu) = \bfU_{\init}(\mu)
    \end{cases}
\end{align}
with the Hamiltonian gradient given by:
\begin{equation}\label{eq:hamiltonian gradient fom}
    \displaystyle\nabla_\bfU \hami\left(\bfU(t;\mu)\right) = \begin{pmatrix}
         \nabla_{\mathbf{x}} \mathcal{U}(\bfX(t;\mu)) \\
         \bfV(t;\mu)
    \end{pmatrix} = 
    \begin{pmatrix}
         \frac{q}{m} \nabla \Lambda^0\left(\bfX(t;\mu)\right) L^{-1} \left( q\omega\Lambda^0\left(\bfX(t;\mu)\right)^T \mathbbm{1}_N - h \rho_0\mathbbm{1}_{n_x} \right) \\
         \bfV(t;\mu)
    \end{pmatrix}
\end{equation}
and $J_{2N}$ referring to the canonical symplectic matrix 
\[ J_{2N} = \begin{pmatrix} 0_N & I_N \\ -I_N & 0_N \end{pmatrix}, \] 
with $I_N$ and $0_N$, respectively, the identity and null matrices of size $N$. Equations \eqref{eq:hamiltonian fom}-\eqref{eq:hamiltonian gradient fom} will be referred to as the Hamiltonian FOM system.

    \subsection{Time discretization and initialization}\label{sec:temporal integration}

We recall that the flow $\phi_t:\R^{2N}\to\R^{2N}$ of a differential equation is a mapping from the initial state to the state at any time $t$
\[
    \phi_t \left(\bfU_{\init}(\mu)\right) := \bfU(t;\mu).
\]
A key property of Hamiltonian systems, as defined in \Cref{eq:hamiltonian fom}, is that the associated flow is symplectic, meaning that it satisfies the relation
\[ \left(\nabla_\bfU \phi_t \left(\bfU_{\init}(\mu)\right)\right)^T J_{2N}\left(\nabla_\bfU \phi_t \left(\bfU_{\init}(\mu)\right)\right) = J_{2N},\qquad \forall t \in [0,T], \mu \in \Gamma. \]
One consequence is that the Hamiltonian $\hami$ is preserved along the flow
\[
\hami \left( \bfU(t;\mu) \right) = \hami \left( \bfU_{\init}(\mu) \right), 
\qquad \forall t \in (0,T], \mu \in \Gamma.
\]
This is particularly important when considering physical systems. To preserve the symplectic structure at the discrete level, we consider the Störmer-Verlet scheme, which is a symplectic time integrator \cite{hairer2006geometric}. It is second order accurate and is explicit in the case of a separable Hamiltonian, which is the case in the problem under consideration. Indeed, the Hamiltonian \eqref{eq:H vlasov poisson} writes as the sum of a discrete kinetic energy, depending only on $\mathbf{v}$, and a discrete potential energy, depending only on $\mathbf{x}$:
\begin{equation*}
    \hami(\bfU) = \hami^{\text{kin}}(\bfV) + \hami^{\text{pot}}(\bfX).
\end{equation*}
with
\begin{align*}
    \hami^{\text{kin}}(\bfV) &= \frac{1}{2} \|\bfV\|^2,\quad \hami^{\text{pot}}(\bfX) = \mathcal{U}(\bfX).
\end{align*}
Introducing a time step $\Delta t$, and denoting $\bfU^n = (\bfX^n, \bfV^n)$ the numerical solution at time $t^n = n \Delta t$, the Störmer-Verlet scheme reads
\begin{equation}\label{eq: stormer verlet step}
    \begin{split}
    \bfV^{n+\half} &= \bfV^n - \frac{\Delta t}{2} \nabla_{\mathbf{x}} \mathcal{U}(\bfX^n), \\
    \bfX^{n+1} &= \bfX^n + \Delta t \, \bfV^{n+\half} , \\ 
    \bfV^{n+1} &= \bfV^{n+1} - \frac{\Delta t}{2} \nabla_{\mathbf{x}} \mathcal{U}(\bfX^{n+1}),\end{split}
\end{equation}
where the expression of $\nabla_{\mathbf{x}}\mathcal{U}$ is given in \Cref{eq:hamiltonian gradient fom}.

The numerical simulation starts by initializing the particle positions, $\bfX^0 = \bfX_{\init}(\mu)$, and velocities, $\bfV^0 = \bfV_{\init}(\mu)$, based on the initial distribution $f_{\init}(x,v;\mu)$. A common approach is to use inverse sampling, which may require to empirical estimate the inverse cumulative distribution function. This method depends on a random number generator, which introduces noise that can degrade the accuracy of the solution \cite{walsh1981nonrandom, barsamian2018picvert}. To avoid this issue, we instead use a non-random number generator based on a Hammersley sequence \cite{hammersley1964random}, which effectively reduces simulation noise. This method is known as a quiet start.

\section{A Hamiltonian reduction with Proper Symplectic Decompostion prereduction}\label{sec: A Hamiltonian reduction with Proper Symplectic Decompostion prereduction}

Taking into consideration that the number of particles $N$ is generally large, the numerical resolution of the Hamiltonian FOM, given in \Cref{eq:hamiltonian fom}, requires significant computational resources and time. Hence, obtaining solutions for various parameters $\mu\in\Gamma$ and times $t$ can become computationally intractable. As a consequence, we aim at building a reduced order model, much smaller in size, that captures the main dynamics for various times $t$ and parameters $\mu\in\Gamma$ and that is more affordable to compute. This reduced order model must also have a Hamiltonian structure.

First, we define the solution manifold
\[
    \mathcal{M} = \left\{ \bfU(t;\mu) \, \text{with} \, t\in[0,T], \mu\in\Gamma \right\} \subset \R^{2N}
\]
formed by the values taken by the solutions of the ODE \Cref{eq:hamiltonian fom}. The manifold structure results from the Cauchy-Lipschitz (Picard-Lindelöf) theorem with parameters under some regularity assumptions of the Hamiltonian. We assume that $\mathcal{M}$ is well approximated by a trial manifold $\widehat{\mathcal{M}}$ that reads
\[
   \widehat{\mathcal{M}} =  \left\{ \calD \left( \bfu(t;\mu) \right)\text{ with } \bfu(t;\mu) \in \R^{2K} \right\} \subset \R^{2N},
\]
with a decoding operator $\calD: \R^{2K} \to \R^{2N}$. In addition, we consider its pseudo-inverse operator $\calE  :\R^{2N} \to \R^{2K}$, called the encoder, which satisfies the relation 
\[ \calE \circ \calD = \text{Id}_{\R^{2K}}. \]
In other words, we search for a reduced model that is a $2K$-dimensional ODE of solution $\bfu(t;\mu)$. To do so, we have to determine $\mathcal{D}$ and $\mathcal{E}$, we therefore ask for the projection operator $ \calD \circ \calE$ onto $\widehat{\mathcal{M}}$ to be close to the identity on a data set $U \subset \mathcal{M}$:
\begin{equation*}
  \forall \bfU \in U,\quad \calD \circ \calE(\bfU) \approx \bfU.
\end{equation*}
 The data set $U$ is composed of snapshots of the solutions at different times  and various parameters, obtained with time integration; it writes
\begin{equation}\label{eq: U snaphots matrix}
    U = \left\{
  \bfU^0_{\mu_1}, \dots ,\bfU^{n_T}_{\mu_1}, \dots ,\bfU^0_{\mu_P}, \dots, \bfU^{n_T}_{\mu_P}
    \right\} \in \mathcal{M}_{2N, (n_T + 1)P}(\R) ,
\end{equation}
where  $\bfU^k_{\mu_p} \approx \bfU(t_k;\mu_p)$ it the numerical solution at time step $k$ and parameters $\mu_p$, $n_T+1 >0$ is the total number of time steps and $P > 0$ is the number of sampled parameters. In practice, parameters are uniformly sampled across $\Gamma$. We denote this sample $\Gamma^{\train}:=\{\mu_p\}_{p\in\{1,\dots,P\}}$.

In addition, we constrain the reduced variables $\bfu(t;\mu)$ to follow the reduced Hamiltonian dynamics
\begin{align}\label{eq:hamiltonian rom}
    \begin{cases}
    \displaystyle \frac{d}{dt} \bfu(t;\mu) = J_{2K} \nabla_\bfu \bar{\hami}\left(\bfu(t;\mu)\right) , &\text{ in } [0,T],\\
    \bfu(0;\mu) = \mathcal{E}\left(\bfU_{\init}(\mu)\right),
    \end{cases}
\end{align}
where $\bar{\hami} : \R^{2K} \to \R$ is a reduced Hamiltonian, to be built.

In the following, we present our strategy to construct the encoder and decoder. It is based on the coupling of the Proper Symplectic Decomposition (PSD), introduced in \cite{peng2015symplectic}, and the AE-HNN method proposed in  \cite{cote2025hamiltonian}, which combines an AutoEncoder (AE) \cite{goodfellow2016deep} and a Hamiltonian Neural Network (HNN) \cite{greydanus2019hamiltonian}. The method will be referred to as PSD-AE-HNN. 

\subsection{PSD-AE-HNN reduction method}\label{sec: PSD and AE-HNN coupled reduced order model}

The PSD-AE-HNN is a three-step reduction method. 

First, the Hamiltonian FOM, which evolves in a $2N$-dimensional phase space, is projected onto an intermediate $2M$-dimensional symplectic subspace with $M \ll N$ using the PSD. Let $A \in \mathcal M_{2N, 2M}(\R)$ denote the symplectic matrix obtained from the PSD algorithm, and $A^+ \in \mathcal M_{2M,2N}(\R)$ be its symplectic inverse, so that $A^+A=I_{2M}$. Together, $A$ and $A^+$ serve as projection and reconstruction operators between the full and intermediate reduced phase space. Further details are provided in \Cref{sec: psd details}.

Second, we further reduce the intermediate $2M$-dimensional representation to a low-dimensional $2K$-dimensional subspace, with $K \ll M$, using an AE. This neural network consists of an encoder $\calE_{\theta_e} : \R^{2M} \to \R^{2K}$ and a decoder $\calD_{\theta_d} : \R^{2K} \to \R^{2M}$, where $\theta_e$ and $\theta_d$ denote the respective parameters of the encoder and decoder. Additional details on the network architectures and training setup are provided in \Cref{sec:Training hyper-parameters}. The autoencoder is trained to approximate the identity mapping, i.e. $\calD_{\theta_d} \circ \calE_{\theta_e} \approx \operatorname{id}$. These networks thus act as nonlinear projectors, mapping data from the intermediate subspace to the low-dimensional space and back. Consequently, the full encoder and decoder are defined by
\[
    \calE = \calE_{\theta_e} \circ A^+, \quad \calD = A \circ \calD_{\theta_d},
\]
where $A^+$ (resp. $A$) is identified with the map $\bfU \mapsto A^+\bfU$ (resp. $\bfU \mapsto A\bfU$).

Third, the dynamics of the reduced model is then captured by a third neural network, the HNN, denoted $\bar{\hami}_{\theta_h}:\R^{2K}\to\R$, where $\theta_h$ represents its trainable parameters. It is trained such that \Cref{eq:hamiltonian rom} holds when evaluated on the reduced variables:
\begin{align}\tag{\ref{eq:hamiltonian rom}}
    \begin{cases}
    \displaystyle \frac{d}{dt} \calE(\bfU(t;\mu)) \approx J_{2K} \nabla_\bfu \bar{\hami}_{\theta_h}\left(\calE(\bfU(t;\mu)) \right) , &\text{ in } [0,T]\\
    \bfu(0;\mu) = \mathcal{E}\left(\bfU_{\init}(\mu)\right)
    \end{cases}
\end{align}
A more detailed description of this component is provided in \Cref{sec: ae-hnn details}.

\begin{figure}
    \centering
\resizebox{1.\textwidth}{!}{

\tikzset{every picture/.style={line width=0.75pt}} 

\begin{tikzpicture}[x=0.75pt,y=0.75pt,yscale=-1,xscale=1]

\draw  [fill={rgb, 255:red, 209; green, 255; blue, 156 }  ,fill opacity=1 ] (200.6,59.65) -- (210.6,59.65) -- (210.6,133.35) -- (200.6,133.35) -- cycle ;
\draw  [fill={rgb, 255:red, 209; green, 255; blue, 156 }  ,fill opacity=1 ] (197,63.65) -- (207,63.65) -- (207,137.35) -- (197,137.35) -- cycle ;

\draw  [fill={rgb, 255:red, 231; green, 130; blue, 255 }  ,fill opacity=1 ] (25.15,23) -- (35.15,23) -- (35.15,174) -- (25.15,174) -- cycle ;
\draw  [fill={rgb, 255:red, 255; green, 179; blue, 189 }  ,fill opacity=1 ] (77.95,42.1) -- (87.95,42.1) -- (87.95,154.9) -- (77.95,154.9) -- cycle ;
\draw  [fill={rgb, 255:red, 209; green, 255; blue, 156 }  ,fill opacity=1 ] (160.71,42.1) -- (170.71,42.1) -- (170.71,154.9) -- (160.71,154.9) -- cycle ;
\draw    (91.2,98.2) -- (151.4,98.49) ;
\draw [shift={(154.4,98.5)}, rotate = 180.27] [fill={rgb, 255:red, 0; green, 0; blue, 0 }  ][line width=0.08]  [draw opacity=0] (8.93,-4.29) -- (0,0) -- (8.93,4.29) -- cycle    ;
\draw    (40.5,98.5) -- (72.2,98.23) ;
\draw [shift={(75.2,98.2)}, rotate = 179.5] [fill={rgb, 255:red, 0; green, 0; blue, 0 }  ][line width=0.08]  [draw opacity=0] (8.93,-4.29) -- (0,0) -- (8.93,4.29) -- cycle    ;
\draw  [fill={rgb, 255:red, 209; green, 255; blue, 156 }  ,fill opacity=1 ] (240.35,77.65) -- (250.35,77.65) -- (250.35,106.15) -- (240.35,106.15) -- cycle ;
\draw  [fill={rgb, 255:red, 209; green, 255; blue, 156 }  ,fill opacity=1 ] (236.75,82.05) -- (246.75,82.05) -- (246.75,110.55) -- (236.75,110.55) -- cycle ;
\draw  [fill={rgb, 255:red, 209; green, 255; blue, 156 }  ,fill opacity=1 ] (233.55,86.45) -- (243.55,86.45) -- (243.55,114.95) -- (233.55,114.95) -- cycle ;
\draw  [fill={rgb, 255:red, 209; green, 255; blue, 156 }  ,fill opacity=1 ] (230.35,90.85) -- (240.35,90.85) -- (240.35,119.35) -- (230.35,119.35) -- cycle ;

\draw  [fill={rgb, 255:red, 209; green, 255; blue, 156 }  ,fill opacity=1 ] (295.15,42.1) -- (305.15,42.1) -- (305.15,154.9) -- (295.15,154.9) -- cycle ;
\draw  [fill={rgb, 255:red, 209; green, 255; blue, 156 }  ,fill opacity=1 ] (311.65,59.25) -- (321.65,59.25) -- (321.65,137.75) -- (311.65,137.75) -- cycle ;
\draw  [fill={rgb, 255:red, 209; green, 255; blue, 156 }  ,fill opacity=1 ] (328.35,75.05) -- (338.35,75.05) -- (338.35,121.95) -- (328.35,121.95) -- cycle ;
\draw  [fill={rgb, 255:red, 209; green, 255; blue, 156 }  ,fill opacity=1 ] (343.55,85.85) -- (353.55,85.85) -- (353.55,111.15) -- (343.55,111.15) -- cycle ;
\draw  [fill={rgb, 255:red, 231; green, 130; blue, 255 }  ,fill opacity=1 ] (586.05,324.1) -- (576.05,324.1) -- (576.05,173.1) -- (586.05,173.1) -- cycle ;
\draw  [fill={rgb, 255:red, 255; green, 179; blue, 189 }  ,fill opacity=1 ] (525.05,304.5) -- (515.05,304.5) -- (515.05,191.7) -- (525.05,191.7) -- cycle ;
\draw  [fill={rgb, 255:red, 187; green, 217; blue, 255 }  ,fill opacity=1 ] (449.56,304.5) -- (439.56,304.5) -- (439.56,191.7) -- (449.56,191.7) -- cycle ;
\draw    (506.4,248.1) -- (455.8,248.1) ;
\draw [shift={(509.4,248.1)}, rotate = 180] [fill={rgb, 255:red, 0; green, 0; blue, 0 }  ][line width=0.08]  [draw opacity=0] (8.93,-4.29) -- (0,0) -- (8.93,4.29) -- cycle    ;
\draw    (568,248.01) -- (531,248.1) ;
\draw [shift={(571,248)}, rotate = 179.86] [fill={rgb, 255:red, 0; green, 0; blue, 0 }  ][line width=0.08]  [draw opacity=0] (8.93,-4.29) -- (0,0) -- (8.93,4.29) -- cycle    ;
\draw  [fill={rgb, 255:red, 187; green, 217; blue, 255 }  ,fill opacity=1 ] (313.05,304.5) -- (303.05,304.5) -- (303.05,191.7) -- (313.05,191.7) -- cycle ;
\draw  [fill={rgb, 255:red, 187; green, 217; blue, 255 }  ,fill opacity=1 ] (295.85,287.35) -- (285.85,287.35) -- (285.85,208.85) -- (295.85,208.85) -- cycle ;
\draw  [fill={rgb, 255:red, 187; green, 217; blue, 255 }  ,fill opacity=1 ] (279.05,271.55) -- (269.05,271.55) -- (269.05,224.65) -- (279.05,224.65) -- cycle ;
\draw  [fill={rgb, 255:red, 187; green, 217; blue, 255 }  ,fill opacity=1 ] (262.65,260.75) -- (252.65,260.75) -- (252.65,235.45) -- (262.65,235.45) -- cycle ;
\draw  [fill={rgb, 255:red, 187; green, 217; blue, 255 }  ,fill opacity=1 ] (373.85,227.25) -- (383.85,227.25) -- (383.85,255.75) -- (373.85,255.75) -- cycle ;
\draw  [fill={rgb, 255:red, 187; green, 217; blue, 255 }  ,fill opacity=1 ] (370.25,231.65) -- (380.25,231.65) -- (380.25,260.15) -- (370.25,260.15) -- cycle ;
\draw  [fill={rgb, 255:red, 187; green, 217; blue, 255 }  ,fill opacity=1 ] (367.05,236.05) -- (377.05,236.05) -- (377.05,264.55) -- (367.05,264.55) -- cycle ;
\draw  [fill={rgb, 255:red, 187; green, 217; blue, 255 }  ,fill opacity=1 ] (363.85,240.45) -- (373.85,240.45) -- (373.85,268.95) -- (363.85,268.95) -- cycle ;

\draw  [fill={rgb, 255:red, 187; green, 217; blue, 255 }  ,fill opacity=1 ] (405.45,209.5) -- (415.45,209.5) -- (415.45,283.2) -- (405.45,283.2) -- cycle ;
\draw  [fill={rgb, 255:red, 187; green, 217; blue, 255 }  ,fill opacity=1 ] (401.85,213.5) -- (411.85,213.5) -- (411.85,287.2) -- (401.85,287.2) -- cycle ;

\draw  [fill={rgb, 255:red, 255; green, 196; blue, 108 }  ,fill opacity=1 ] (72.95,235.45) -- (82.95,235.45) -- (82.95,260.75) -- (72.95,260.75) -- cycle ;
\draw  [fill={rgb, 255:red, 255; green, 196; blue, 108 }  ,fill opacity=1 ] (91.55,224.65) -- (101.55,224.65) -- (101.55,271.55) -- (91.55,271.55) -- cycle ;
\draw  [fill={rgb, 255:red, 255; green, 196; blue, 108 }  ,fill opacity=1 ] (109.05,224.65) -- (119.05,224.65) -- (119.05,271.55) -- (109.05,271.55) -- cycle ;
\draw  [fill={rgb, 255:red, 255; green, 196; blue, 108 }  ,fill opacity=1 ] (127.05,224.65) -- (137.05,224.65) -- (137.05,271.55) -- (127.05,271.55) -- cycle ;
\draw  [fill={rgb, 255:red, 255; green, 196; blue, 108 }  ,fill opacity=1 ] (144.55,241.38) -- (154.55,241.38) -- (154.55,254.83) -- (144.55,254.83) -- cycle ;
\draw    (246.5,247.81) -- (158.5,248) ;
\draw [shift={(249.5,247.8)}, rotate = 179.87] [fill={rgb, 255:red, 0; green, 0; blue, 0 }  ][line width=0.08]  [draw opacity=0] (8.93,-4.29) -- (0,0) -- (8.93,4.29) -- cycle    ;
\draw [line width=0.75]    (358.5,97.5) .. controls (379.5,97.5) and (383.5,139.5) .. (360,154) .. controls (336.5,168.5) and (141.5,162.5) .. (110.5,172) .. controls (79.5,181.5) and (60.5,194.5) .. (60,212.5) .. controls (59.53,229.51) and (55.48,238.48) .. (68.12,248.28) ;
\draw [shift={(70.5,250)}, rotate = 214.11] [fill={rgb, 255:red, 0; green, 0; blue, 0 }  ][line width=0.08]  [draw opacity=0] (8.93,-4.29) -- (0,0) -- (8.93,4.29) -- cycle    ;
\draw  [fill={rgb, 255:red, 65; green, 117; blue, 5 }  ,fill opacity=1 ][line width=0.75]  (160.71,63.67) -- (170.71,63.67) -- (170.71,100.33) -- (160.71,100.33) -- cycle ;
\draw  [fill={rgb, 255:red, 65; green, 117; blue, 5 }  ,fill opacity=1 ][line width=0.75]  (196.92,98.67) -- (206.92,98.67) -- (206.92,110.33) -- (196.92,110.33) -- cycle ;
\draw  [fill={rgb, 255:red, 7; green, 71; blue, 151 }  ,fill opacity=1 ][line width=0.75]  (401.67,221.25) -- (411.67,221.25) -- (411.67,232.92) -- (401.67,232.92) -- cycle ;
\draw  [fill={rgb, 255:red, 7; green, 71; blue, 151 }  ,fill opacity=1 ][line width=0.75]  (439.56,244) -- (449.56,244) -- (449.56,277) -- (439.56,277) -- cycle ;
\draw [color={rgb, 255:red, 70; green, 134; blue, 0 }  ,draw opacity=1 ][fill={rgb, 255:red, 65; green, 117; blue, 5 }  ,fill opacity=1 ]   (170.71,63.67) -- (196.92,98.67) ;
\draw [color={rgb, 255:red, 70; green, 134; blue, 0 }  ,draw opacity=1 ][fill={rgb, 255:red, 65; green, 117; blue, 5 }  ,fill opacity=1 ]   (170.71,100.33) -- (196.92,110.33) ;
\draw [color={rgb, 255:red, 20; green, 99; blue, 197 }  ,draw opacity=1 ][fill={rgb, 255:red, 20; green, 99; blue, 197 }  ,fill opacity=1 ]   (411.67,221.25) -- (439.56,244) ;
\draw [color={rgb, 255:red, 20; green, 99; blue, 197 }  ,draw opacity=1 ][fill={rgb, 255:red, 20; green, 99; blue, 197 }  ,fill opacity=1 ]   (439.56,277) -- (411.67,232.92) ;
\draw  [color={rgb, 255:red, 208; green, 2; blue, 27 }  ,draw opacity=1 ] (53,25.5) .. controls (53,17.77) and (59.27,11.5) .. (67,11.5) -- (532,11.5) .. controls (539.73,11.5) and (546,17.77) .. (546,25.5) -- (546,326) .. controls (546,333.73) and (539.73,340) .. (532,340) -- (67,340) .. controls (59.27,340) and (53,333.73) .. (53,326) -- cycle ;
\draw    (254.8,99.5) -- (288.6,99.78) ;
\draw [shift={(291.6,99.8)}, rotate = 180.47] [fill={rgb, 255:red, 0; green, 0; blue, 0 }  ][line width=0.08]  [draw opacity=0] (8.93,-4.29) -- (0,0) -- (8.93,4.29) -- cycle    ;
\draw    (357,247.49) -- (316.4,247.4) ;
\draw [shift={(360,247.5)}, rotate = 180.13] [fill={rgb, 255:red, 0; green, 0; blue, 0 }  ][line width=0.08]  [draw opacity=0] (8.93,-4.29) -- (0,0) -- (8.93,4.29) -- cycle    ;
\draw  [fill={rgb, 255:red, 126; green, 211; blue, 33 }  ,fill opacity=1 ][line width=0.75]  (196.92,110.33) -- (206.92,110.33) -- (206.92,127.17) -- (196.92,127.17) -- cycle ;
\draw  [fill={rgb, 255:red, 126; green, 211; blue, 33 }  ,fill opacity=1 ][line width=0.75]  (230.38,101.5) -- (240.17,101.5) -- (240.17,112.92) -- (230.38,112.92) -- cycle ;
\draw [color={rgb, 255:red, 126; green, 211; blue, 33 }  ,draw opacity=1 ][fill={rgb, 255:red, 126; green, 211; blue, 33 }  ,fill opacity=1 ]   (206.92,127.17) -- (230.38,112.92) ;
\draw [color={rgb, 255:red, 126; green, 211; blue, 33 }  ,draw opacity=1 ][fill={rgb, 255:red, 126; green, 211; blue, 33 }  ,fill opacity=1 ]   (206.92,98.67) -- (230.38,101.5) ;
\draw  [fill={rgb, 255:red, 80; green, 227; blue, 194 }  ,fill opacity=1 ][line width=0.75]  (363.88,250.75) -- (373.67,250.75) -- (373.67,262.42) -- (363.88,262.42) -- cycle ;
\draw  [fill={rgb, 255:red, 80; green, 227; blue, 194 }  ,fill opacity=1 ][line width=0.75]  (401.88,241.25) -- (411.67,241.25) -- (411.67,272.75) -- (401.88,272.75) -- cycle ;
\draw [color={rgb, 255:red, 80; green, 227; blue, 194 }  ,draw opacity=1 ][fill={rgb, 255:red, 80; green, 227; blue, 194 }  ,fill opacity=1 ]   (401.88,241.25) -- (373.67,250.75) ;
\draw [color={rgb, 255:red, 80; green, 227; blue, 194 }  ,draw opacity=1 ][fill={rgb, 255:red, 80; green, 227; blue, 194 }  ,fill opacity=1 ]   (401.88,272.75) -- (373.67,262.42) ;

\draw (5,7) node [anchor=north west][inner sep=0.75pt]  [font=\small]  {$\mathbf{u}( t_{1} ;\mu )$};
\draw (56.6,25) node [anchor=north west][inner sep=0.75pt]  [font=\small]  {$A^{+}\mathbf{u}( t_{1} ;\mu )$};
\draw (343,70) node [anchor=north west][inner sep=0.75pt]  [font=\small]  {$\bar{\mathbf{u}}( t_{1} ;\mu )$};
\draw (222,219) node [anchor=north west][inner sep=0.75pt]  [font=\small]  {$\bar{\mathbf{u}}( t_{2} ;\mu )$};
\draw (494.3,175) node [anchor=north west][inner sep=0.75pt]  [font=\small]  {$\tilde{\mathbf{u}}( t_{2} ;\mu )$};
\draw (559.2,157) node [anchor=north west][inner sep=0.75pt]  [font=\small]  {$A\tilde{\mathbf{u}}( t_{2} ;\mu )$};
\draw (91,275.5) node [anchor=north west][inner sep=0.75pt]  [font=\small,color={rgb, 255:red, 255; green, 162; blue, 0 }  ,opacity=1 ]  {$\text{HNN} \ \hami_{\theta _{h}}$};
\draw (91,83) node [anchor=north west][inner sep=0.75pt]  [font=\scriptsize] [align=left] {expand dim.};
\draw (459,238.5) node [anchor=north west][inner sep=0.75pt]  [font=\scriptsize] [align=left] {squeeze};
\draw (162.5,250) node [anchor=north west][inner sep=0.75pt]  [font=\scriptsize] [align=left] {time integration};
\draw (252,304) node [anchor=north west][inner sep=0.75pt]  [font=\small,color={rgb, 255:red, 0; green, 106; blue, 231 }  ,opacity=1 ]  {$\underbrace{\ \ \ \ \ \ \ \ \ \ \ \ \ \ \ \ \ \ \ \ \ \ \ \ \ \ \ \ \ \ \ \ \ \ \ \ \ \ \ \ \ \ \ \ \ \ \ \ \ \ }$};
\draw (355,43) node [anchor=north west][inner sep=0.75pt]  [font=\small,color={rgb, 255:red, 102; green, 201; blue, 0 }  ,opacity=1 ,rotate=-180]  {$\underbrace{\ \ \ \ \ \ \ \ \ \ \ \ \ \ \ \ \ \ \ \ \ \ \ \ \ \ \ \ \ \ \ \ \ \ \ \ \ \ \ \ \ \ \ \ \ \ \ \ \ \ }$};
\draw (221,20) node [anchor=north west][inner sep=0.75pt]  [font=\small,color={rgb, 255:red, 102; green, 201; blue, 0 }  ,opacity=1 ]  {$\text{encoder} \ \mathcal{E}_{\theta _{e}}$};
\draw (319,316) node [anchor=north west][inner sep=0.75pt]  [font=\small,color={rgb, 255:red, 0; green, 106; blue, 231 }  ,opacity=1 ]  {$\text{decoder} \ \mathcal{D}_{\theta _{d}}$};
\draw (121,343) node [anchor=north west][inner sep=0.75pt]  [font=\small,color={rgb, 255:red, 208; green, 2; blue, 27 }  ,opacity=1 ] [align=left] {PSD-AE-HNN method};
\draw (253,85) node [anchor=north west][inner sep=0.75pt]  [font=\scriptsize] [align=left] {flatten};
\draw (313,233) node [anchor=north west][inner sep=0.75pt]  [font=\scriptsize] [align=left] {unflatten};

\end{tikzpicture}
}    
    \caption{PSD-AE-HNN architecture: from FOM solution $\bfU(t_1;\mu)$, a PSD intermediate reduced variable $A^+ \bfU(t_1;\mu)$ is computed, followed by the reduced state $\bar{\mathbf{u}}( t_{1} ;\mu )=\calE_{\theta_e}(A^+ \bfU(t_1;\mu))$. Next,  time integration $\bar{\mathbf{u}}( t_{2} ;\mu )$ is performed with the HNN gradient, the final state $A\tilde{\mathbf{u}}( t_{2} ;\mu )=A\calD_{\theta_d}(\bar{\mathbf{u}}( t_{2} ;\mu ))$ is recovered with decoder and PSD successive decompression.}\label{fig:scheme PSD+AE-HNN}
\end{figure}
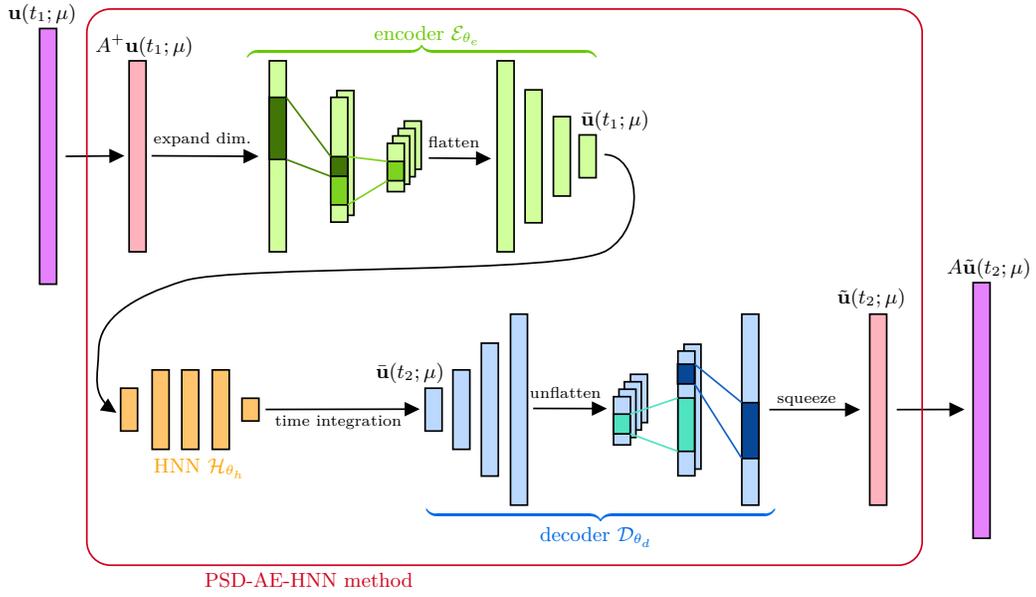

The online process for applying the reduced model is schematized in \Cref{fig:scheme PSD+AE-HNN}. We start with a full order solution $\bfU(t_1;\mu) \in \mathbb{R}^{2N}$ at time $t_1$. Our goal is to approximate the full order solution at time $t_2 > t_1$, using the reduced model. We first apply the symplectic projection to an intermediate reduced variable 
\[
A^+ \bfU(t_1;\mu) \in \mathbb{R}^{2M}.
\]
The encoder then maps this intermediate representation to a low-dimensional reduced state
\[
\bfu(t_1;\mu) =\mathcal{E}_{\theta_e}(A^+ \bfU(t_1;\mu))\in \mathbb{R}^{2K}.
\]
Since $\bfu(t;\mu)$ evolves according to a Hamiltonian system defined by the HNN, we employ the Störmer-Verlet integrator described in \Cref{eq: stormer verlet step} to advance the solution in time up to time $t_2$. The required gradients of the learned Hamiltonian are computed via backpropagation, allowing us to obtain the reduced state $\bfu(t_2;\mu)$ at time $t_2$.
Finally, the decompression step is performed to recover an approximation of the full-order solution. The reduced state $\bfu(t_2;\mu)$ is first decoded to the intermediate space via
\[
\tilde{\mathbf{u}}( t_{2} ;\mu )=\mathcal{D}_{\theta_d}(\bfu(t_2;\mu)).
\]
Finally, we apply the symplectic lift to reconstruct the full-order approximation
\[
A\tilde{\mathbf{u}}( t_{2} ;\mu ) \approx \bfU ( t_{2} ;\mu ).
\]

There are two main motivations for combining the PSD with the AE-HNN. First, although the AE-HNN is an efficient data-driven model reduction technique, its computational cost scales with its input dimension. For large $N$, this results in neural networks that are too large to train effectively. Second, since the inputs correspond to particles in phase space, they are inherently unstructured and may contain noise. The prior reduction via PSD projects the dynamics onto a lower-dimensional symplectic subspace, resulting in a more structured and compact representation. This intermediate reduced variable is both easier to learn for the autoencoder and HNN while also preserving the underlying Hamiltonian structure.

The offline stage of the method, to construct the different elements of the reduced order model, consists of three main steps:

(i) snapshot generation: we compute a collection of full order solutions at various times and for different parameter values;

(ii) symplectic basis construction: we apply the PSD algorithm to build the reduced symplectic basis $A$;

(iii) neural network training: following the approach of \cite{cote2025hamiltonian}, we simultaneously train the second stage of the encoder, $\calE_{\theta_e}$, the first stage of the decoder, $\calD_{\theta_d}$, and the HNN, $\bar{\hami}_{\theta_h}$. These networks are trained using the FOM snapshots projected onto the intermediate subspace via $A^+$.

We dive into both PSD and AE-HNN functioning in the following sections.    
    
    \subsection{PSD reduction}\label{sec: psd details}

In this section, we briefly introduce the Proper Symplectic Decomposition (PSD) \cite{peng2015symplectic}. The goal of PSD is to approximate the manifold $\mathcal{M} \subset \R^{2N}$ of full order states with a $2M$-dimensional linear subspace. To preserve the Hamiltonian structure of the dynamics, we require the projection to be symplectic. Under this constraint, the intermediate reduced variable $\widetilde{\bfU}(t;\mu) \in \R^{2M}$ is defined by
\begin{equation} \label{def:tilde_bfU_2}
\widetilde{\bfU}(t;\mu) = A^+ \bfU(t;\mu),
\end{equation}
where $A^+$ denotes the symplectic inverse of a matrix $A \in \operatorname{Sp}_{2M, 2N}(\R)$. This set denotes the symplectic Stiefel manifold, which consists of all $2N \times 2M$ matrices $A$ satisfying the symplectic condition
\[
    A^T J_{2N} A = J_{2M}.
\]
For any matrix $A \in \operatorname{Sp}_{2M, 2N}(\R)$, its symplectic inverse $A^+$ is given by
\begin{equation} \label{def:A+}
    A^+ = J_{2M}^T A^T J_{2N},
\end{equation}
which satisfies $A^+A=I_{2M}$. 

The symplectic matrix $A$ is computed by minimizing the reconstruction error over a set of training snapshots. That is, $A$ is obtained as the solution to the following optimization problem
\begin{equation}\label{eq: min problem PSD 1}
    \underset{A \in \operatorname{Sp}_{2M, 2N}(\R)}{\operatorname{min}} \left\| U - A A^+ U \right\|_F,
\end{equation}
where $\|X\|_F := \sqrt{ \sum_{i,j} | x_{i,j} |^2 }$ is the Frobenius norm and $U$ is the snapshot matrix defined in \Cref{eq: U snaphots matrix}. A direct solution of \Cref{eq: min problem PSD 1} cannot be obtained. However, with additional assumptions outlined in \Cref{annex : PSD details}, we construct the matrix $A$ using the Complex Singular Value Decomposition (SVD) algorithm \cite{peng2015symplectic}. Since $A$ is a symplectic transformation, it can be checked that the intermediate reduced variable $\widetilde{\bfU}$ evolves according to the Hamiltonian dynamics with Hamiltonian function $\hami \circ A$:
\begin{equation}\label{eq: psd reduced model}
\begin{cases}
    \displaystyle \frac{d}{dt} \widetilde{\bfU}(t;\mu) = J_{2M} \nabla_{\tilde{\bfU}} (\hami \circ A)\left(\widetilde{\bfU}(t;\mu)\right) =  J_{2M} A^T \nabla_\bfU \hami\left(A \bfu(t;\mu)\right), \\
    \widetilde{\bfU}(0;\mu) = A^+ \bfU_{\init}(\mu).
\end{cases}
\end{equation}
As observed in \Cref{sec: numerical results}, the value of $M$ required to achieve satisfactory precision is often too large, which reduces the efficiency of a reduced model based solely on the PSD. 
Additionally, the evaluation of the reduced Hamiltonian gradients, $A^T \nabla_\bfU \hami(A \, \cdot)$, still depends on the gradient of the original Hamiltonian function. This results in a computational cost that is higher than that of the FOM itself. To address this issue, hyper-reduction techniques have been proposed, as in \cite{hesthaven2024adaptive}.

    \subsection{AE-HNN reduction}\label{sec: ae-hnn details}

This section provides a brief overview of the AE-HNN method introduced in \cite{cote2025hamiltonian}. The method consists of training simultaneously an auto-encoder, composed of $\calE_{\theta_e}, \calD_{\theta_d}$, and a Hamiltonian Neural Network, $\bar{\hami}_{\theta_h}$. The AE consists of a pair of convolutional neural networks, with convolutional layers followed by dense layers, while the HNN is implemented as a dense neural network \cite{kubat1999neural, goodfellow2016deep}. The neural network parameters, $(\theta_e,\theta_d,\theta_h)  \in \Theta$, are determined by solving an optimization problem of the form
\[
    \underset{(\theta_e,\theta_d,\theta_h)  \in \Theta}{\operatorname{argmin}}\quad  \mathcal{L}(\theta_e,\theta_d,\theta_h),
\]
where the loss function $\mathcal{L}$ is computed using the training dataset $\widetilde{U}$, composed of the snapshots $U$ projected onto the intermediate subspace
\[
   \widetilde{U} = \left\{
  \widetilde{\bfU}^0_{\mu_1}, \dots, \widetilde{\bfU}^{n_T}_{\mu_P}
    \right\}
    = \left\{
  A^+ \bfU^0_{\mu_1}, \dots, A^+ \bfU^{n_T}_{\mu_P}
    \right\}.
\] A gradient descent algorithm is used to determine optimal parameters.

In the AE-HNN method, the loss function is composed of four different loss terms. 
The first term, $\mathcal{L}_{\AE}$, forces the AE to be close to the identity map, i.e.  $\calD_{\theta_d} \circ \calE_{\theta_e} \approx \operatorname{id}$, on the training dataset:
\begin{equation}\label{eq: AE_loss}
    \mathcal{L}_{\AE}(\theta_e, \theta_d) = \sum_{\widetilde{\bfU} \,\in\, \widetilde{U}} \left\| \widetilde{\bfU} - (\calD_{\theta_d}\circ \calE_{\theta_e}) \left( \widetilde{\bfU} \right) \right\|^2_2.
\end{equation}
In practice, a split AE is employed where both the encoder and decoder are made of two neural networks. The first network processes the generalized positions, while the second network processes the generalized velocities. Specifically, the encoder and decoder are structured as follows
\[
    \calE_{\theta_e} = \begin{pmatrix}
    \calE_{\theta_{e,1}} \\ \calE_{\theta_{e,2}} \end{pmatrix}, \quad
    \calD_{\theta_d} = \begin{pmatrix}
    \calD_{\theta_{d,1}} \\ \calD_{\theta_{d,2}} \end{pmatrix}.
\]
The reduced state is then given by
\[
    \bfu(t;\mu)
    = \begin{pmatrix}
        \bfx(t;\mu) \\ \bfv(t;\mu)
    \end{pmatrix}
    = \begin{pmatrix}
        \calE_{\theta_{e,1}}(\widetilde{\bfX}(t;\mu)) \\ \calE_{\theta_{e,2}}(\widetilde{\bfV}(t;\mu))
    \end{pmatrix}
\]
and conversely for the decoded state. 

The second loss term is defined to constrain the reduced trajectories $\bfu(t;\mu)$ to be close to those of a Hamiltonian system, as described in \Cref{eq:hamiltonian rom}. In practice, these reduced dynamics are defined through a time discretization. We therefore introduce the prediction operator
\[
    \mathcal{P}_s\left(\bfu, \bar{\hami}_{\theta_h} \right)
\]
which consists in performing $s \in \mathbb{N}^\ast$ iterations of the Störmer-Verlet algorithm \eqref{eq: stormer verlet step}, starting from the reduced state $\bfu$ and using the reduced Hamiltonian $\bar{\hami}_{\theta_h}$. The number of  iterations considered $s$, also called the watch duration, is a hyperparameter, which must be set. With this prediction operator, the second loss function,  $\mathcal{L}^s_{\modelr}$, constrains the HNN to accurately capture the reduced dynamics between the $n$-th and $(n+s)$-th time steps
\begin{equation}\label{eq: HNN_loss}
  \mathcal{L}^s_{\modelr} (\theta_e,\theta_h) = \sum_{ \widetilde{\bfU}^{n},  \widetilde{\bfU}^{n+s} \in \, \widetilde{U}}\left\| \bfu^{n+s} - \mathcal{P}_s\left( \bfu^{n} ; \bar\hami_{\theta_h} \right) \right\|^2_2,
\end{equation}
where $\widetilde{\bfU}^{n},  \widetilde{\bfU}^{n+s} \in \, \widetilde{U}$ denotes the sampling of random pairs on the dataset $\widetilde{U}$.
Since the full order Hamiltonian in \Cref{eq:H vlasov poisson} is separable, the reduced Hamiltonian $\bar{\hami}_{\theta_h}$ is also assumed to be separable: 
\[
\bar{\hami}_{\theta_h}(\bfu) = \bar{\hami}_{\theta_{h}}^{\text{kin}}(\bfv) + \bar{\hami}_{\theta_{h}}^{\text{pot}} (\bfx),
\]
which allows for the explicit formulation of the time integrator $\mathcal{P}_s$. 

The third part of the loss function aims at ensuring that the reduced trajectories, generated by the encoder $\calE_{\theta_e}$,  preserve the reduced Hamiltonian. The loss function, $\mathcal{L}^s_{\stab}$ writes:
\begin{equation}\label{eq: energy_loss}
    \mathcal{L}^s_{\stab} (\theta_e, \theta_h) = \sum_{ \bar{\bfU}^{n},  \bar{\bfU}^{n+s} \in \, \bar{U}} \left\| \bar\hami_{\theta_h} \left( \bar{\bfU}^{n+s} \right) - \bar\hami_{\theta_h} \left( \bar{\bfU}^{n} \right) \right\|^2_2.
\end{equation}

Finally, the three neural networks are strongly coupled using a loss function, $\mathcal{L}_{\model}^s$, which encapsulates the full prediction from the $n$-th time step to the $n+s$-th time step, using the encoder at the beginning, the decoder at the end and the prediction operator associated with the reduced model:
\begin{equation}\label{eq: model_loss}
\mathcal{L}_{\model}^s (\theta_e,\theta_d,\theta_h) = \sum_{\widetilde{\bfU}^{n},  \widetilde{\bfU}^{n+s} \in \, \widetilde{U}} \left\| \widetilde{\bfU}^{n+s} - \mathcal{D}_{\theta_d} \left( \mathcal{P}_s\left( \mathcal{E}_{\theta_e}(\widetilde{\bfU}^n); \bar\hami_{\theta_h} \right) \right) \right\|^2_2.
\end{equation}

To summarize, four different loss functions, given in \Cref{eq: AE_loss,eq: HNN_loss,eq: energy_loss,eq: model_loss}, are used to train the AE and the HNN. More precisely, the parameters are determined such as to minimize the following weighted sum
\begin{align*}
   \mathcal{L}(\theta_e,\theta_d,\theta_h) = &
    \omega_{\AE}\, \mathcal{L}_{\AE}(\theta_e, \theta_d) 
    + \omega_{\modelr}\, \mathcal{L}^s_{\modelr}(\theta_e, \theta_h)\\
    &+ \omega_{\stab}\, \mathcal{L}^s_{\stab}(\theta_e, \theta_h)
    + \omega_{\model}\, \mathcal{L}_{\model}^s(\theta_e,\theta_d,\theta_h),
\end{align*}
where $\omega_{\AE}, \omega_{\modelr}, \omega_{\stab}$ and $\omega_{\model}$ are positive weights. The networks are thus jointly trained, with potentially adversarial goals. Ultimately, $\mathcal{L}_{\model}^s$ serves as the primary loss function to measure the performance of the AE-HNN reduction. The other loss functions act as auxiliary functions to drive the training process.

    \subsection{Hyperparameters tuning}\label{sec:Training hyper-parameters}

This section specifies the hyperparameters of the models and describes how they are selected. They are chosen based on two main criteria. First, the reduced model must closely approximate the full model, with the difference measured by the losses, while minimizing the reduced dimension $K$. Second, the networks, particularly the HNN, must remain lightweight in terms of the number of parameters to ensure fast computation. Note that the AE is less critical in terms of size, as it is only called once during the online phase.

Regarding the PSD part, the main hyperparameter is the intermediate subspace dimension $M$. A smaller value of $M$ results in a more significant reduction and reduces the computation time, while a larger value of $M$ provides a richer subspace for the subsequent training the AE-HNN, but with increased computation time. In practice, we select an $M$ value such that the PSD reconstruction error is slightly less than the target accuracy of the reduced model. In the following test cases, a typical value is $M=121$ for a final time $T=20$ and $M=256$ when $T=40$.

Secondly, the AE-HNN part involves hyperparameters for defining the architecture of the neural networks. As explained in \Cref{sec: ae-hnn details}, the encoder consists of two convolutional neural networks when considering a split AE. Each network starts with an input of size $M$, which is fed through a series of 1D convolutional layers with a stride of 3, a kernel size of 3, and valid padding each. The number of filters is progressively multiplied by the stride between layers, starting with 12 filters.  The final output is flattened and passed through a series of dense layers, whose sizes gradually decrease, ultimately leading to a single dense layer of output size $K$. The activation function is applied throughout the encoder, except for the output layer, which uses a linear activation function. The decoder is designed as a mirror image of the encoder, where the 1D convolutions are replaced with 1D transposed convolutions. Lastly, the HNN is a simple multi-layer perceptron with an input size of $2K$ and an output size of 1. The activation function in the HNN may differ from that used in the AE.

The AE-HNN also requires some hyperparameters to be fixed for the training. The chosen optimization method is the Adam algorithm \cite{kingma2017adam}, which is an  adaptive stochastic gradient descent method. 
The learning rate follows the rule
\[
    \rho_k =  (0.99)^{k / 150}\, \rho_0,
\]
where the division operator denotes integer division and $k$ is the training step. 
Additionally, we can reset the decay, i.e. set $k=0$, if the loss function reaches a plateau. The purpose of this reset strategy is to escape poor local minima by introducing a sudden, larger learning rate. In most cases, we start training with a large $\rho_0 = 10^{-3}$ to accelerate the convergence. Then, we diminish it to $\rho_0 = 5 \sci{-4}$ or so for fine-tuning. The training process depends on the watch duration $s$. It is be set to $s = 8$ and then be reasonably increased up to $s = 32$ to improve predictions. Finally,training is divided in two stages. First, the AE is trained alone by setting
\[
    \omega_{\AE}=1,\quad \omega_{\modelr}=\omega_{\stab}=\omega_{\model}=0.
\]
Then, after the loss has reached a value in the range $[5\sci{-3}, 1\sci{-2}]$, the AE and the HNN are trained together by setting
\[
    \omega_{\AE}=1,\quad \omega_{\modelr}=10,\quad \omega_{\stab}=1\sci{-4},\quad \omega_{\model}=1.
\]
Table \ref{tab:hyperparameters pic hnn} recapitulates the hyperparameters used for the different test cases of the next section. 

\begin{table}[htb!]
    \centering
    \resizebox{\columnwidth}{!}{%
    \begin{tabular}{cp{4cm}ccc}
    \toprule
    && linear Landau damping & nonlinear Landau damping & two stream instability \\
    \midrule
    AE   & nb of convolution blocks (encoder) & $2$ & $2$ & $2$ \\
         & dense layers (encoder) &  $[150, 100, 50, 25]$ & $[250, 150, 100, 50, 25]$ & $[150, 100, 50, 25]$ \\
         &  activation functions & ELU & ELU & ELU \\
    \midrule
    HNN & dense layers & $[48,24,24,24,12]$ & $[96,48,48,48,24]$ & $[48,24,24,24,12]$ \\
        &  activation functions & softplus & softplus & softplus \\
     \midrule
     watch duration & $s$ & $16$ & $10\to22$ & $16\to32$  \\
    \bottomrule
    \end{tabular}}%
    \caption{Hyper-parameters. Activation functions are used except for the last layer of the neural networks. ELU refers to the function elu$(x) = x 1_{x>0} + (e^x-1) 1_{x<0}$ and softplus to the function softplus$(x) = \log(1+e^{x})$. For the autoencoder (AE), the number of convolution blocks and the sizes of the hidden of layers are those of the encoder. The decoder is constructed in a mirror way.}
    \label{tab:hyperparameters pic hnn}
\end{table}

\begin{remark}
In practice, pre-processing is applied to the neural network inputs. While such functions could be learned by the first layers of the network, manually selecting them significantly improves both the training and prediction processes. Considering the SVD of the snapshot matrix $U$ defined in \Cref{eq: U snaphots matrix},
$$U = W \Sigma V^*,$$ with $W$ and $V$ unitary matrices, $V^*$ is the conjugate transpose of $V$ and $\Sigma$ a diagonal matrix of diagonal values $(\sigma_k)_k$ sorted in descending order, the encoder input is pre-processed with the function
\[
    (\widetilde{\bfU})_k \mapsto \sigma_k^{-1/2}  (\widetilde{\bfU})_k,
\]
where $(\widetilde{\bfU})_k$ is the $k$-th coefficient of the intermediate reduced variable $\widetilde{\bfU}$. The idea is to balance the influence of each singular PSD vector in the intermediate reduced basis, thereby allowing the AE to capture the most important modes without overly neglecting the other modes.
\end{remark}
    
\section{Numerical results}\label{sec: numerical results}

In this section, the PSD-AE-HNN reduction of the PIC method is tested on three classical plasma physics dynamics: the linear Landau damping, the nonlinear Landau damping, and the two-stream instability test cases. 

The parameterized initial distributions of the particles takes the following form
\[
    f_{\init}(x,v;\mu) = f_{\init, x}(x;\alpha) \, f_{\init, v}(v;\sigma),
\]
with parameters $\mu=(\alpha, \sigma)^T \in \Gamma \subset \R^2$. The initial position distribution is a perturbed uniform distribution
\begin{equation} \label{def:f_init_x}
    f_{\init, x}(x;\alpha) = \frac{k}{2\pi} \left(1 + \alpha \cos (k\,x)\right),
\end{equation}
defined over $\Omega_x = \left[0, \frac{2\pi}{k}\right)$, where $k>0$ is a fixed wave number. The parameter $\alpha>0$ is the perturbation amplitude. The initial velocity distribution $f_{\init, v}(v;\sigma)$, defined over $\Omega_v=[-6,6]$, is given by a Gaussian
\begin{equation} \label{def:f_init_v_landaudamping}
    f_{\init,v}(v;\sigma) = \frac{1}{\sigma \sqrt{2\pi}} \operatorname{exp} \left( - \frac{v^2}{2 \sigma^2} \right),
\end{equation}
for the Landau test cases, and by the sum of two Gaussian  
\begin{equation} \label{def:f_init_v_twostream}
    f_{\init,v}(v;\sigma) = \frac{1}{2 \sigma \sqrt{2\pi}} \left[ \operatorname{exp} \left( - \frac{(v - 3)^2}{2 \sigma^2} \right) + \operatorname{exp} \left( - \frac{(v + 3)^2}{2 \sigma^2} \right) \right],
\end{equation}
for the two stream instability test case, where $\sigma > 0$ stands for the standard deviation of the Gaussian distributions.

The $P$ training parameters $\mu\in\Gamma^{\train}$ are selected on a $\sqrt{P} \times \sqrt{P}$ grid over $\Gamma$. The model is then evaluated on a fine $20\times20$ grid $\Gamma^{\test} \subset \Gamma$. For each parameter $\mu$, the reference FOM solution is denoted
\[
    X^{\reference}_\mu =\left\{\bfX^0_{\mu}, \dots ,\bfX^{n_T}_{\mu}\right\}, \quad V^{\reference}_\mu =\left\{\bfV^0_{\mu}, \dots ,\bfV^{n_T}_{\mu}\right\},
\]
while the solution obtained by the PSD-AE-HNN method is denoted
\[
    X^{\test}_\mu =\left\{\hat\bfX^0_{\mu}, \dots ,\hat\bfX^{n_T}_{\mu}\right\}, \quad V^{\test}_\mu =\left\{\hat\bfV^0_{\mu}, \dots ,\hat\bfV^{n_T}_{\mu}\right\}.
\]
We recall that it is obtained through the compression of the initial condition, its complete integration over $[0,T]$ using the HNN followed by its decompression.
We measure the relative errors on a single parameter $\mu$ for all time steps 
\[
    \error_{X, \mu} = \frac{\left\| X^{\test}_\mu  - X^{\reference}_\mu \right\|_F}{\left\|X^{\reference}_\mu \right\|_F},
    \quad
    \error_{V, \mu} = \frac{\left\| V^{\test}_\mu  - V^{\reference}_\mu \right\|_F}{\left\|V^{\reference}_\mu \right\|_F},
\]
and the mean relative errors at a single time $t$ over all $\mu \in \Gamma^{\test}$
\[
    \error^{\opmean}_{X, t} = \opmean \left(\frac{\left\| \bfX^t_\mu  - \hat\bfX^t_\mu \right\|_F}{\left\| \bfX^t_\mu  \right\|_F}, \mu \in \Gamma^{\test} \right),
\quad
    \error^{\opmean}_{V, t} = \opmean \left(\frac{\left\| \bfV^t_\mu  - \hat\bfV^t_\mu \right\|_F}{\left\| \bfV^t_\mu  \right\|_F}, \mu \in \Gamma^{\test} \right).
\]
In addition, we also compute the associated maximal and minimal errors. 

\subsection{Linear Landau damping}\label{sec: linear landau damping}

\begin{figure}[htb!]
    \centering
    \includegraphics[width=\textwidth]{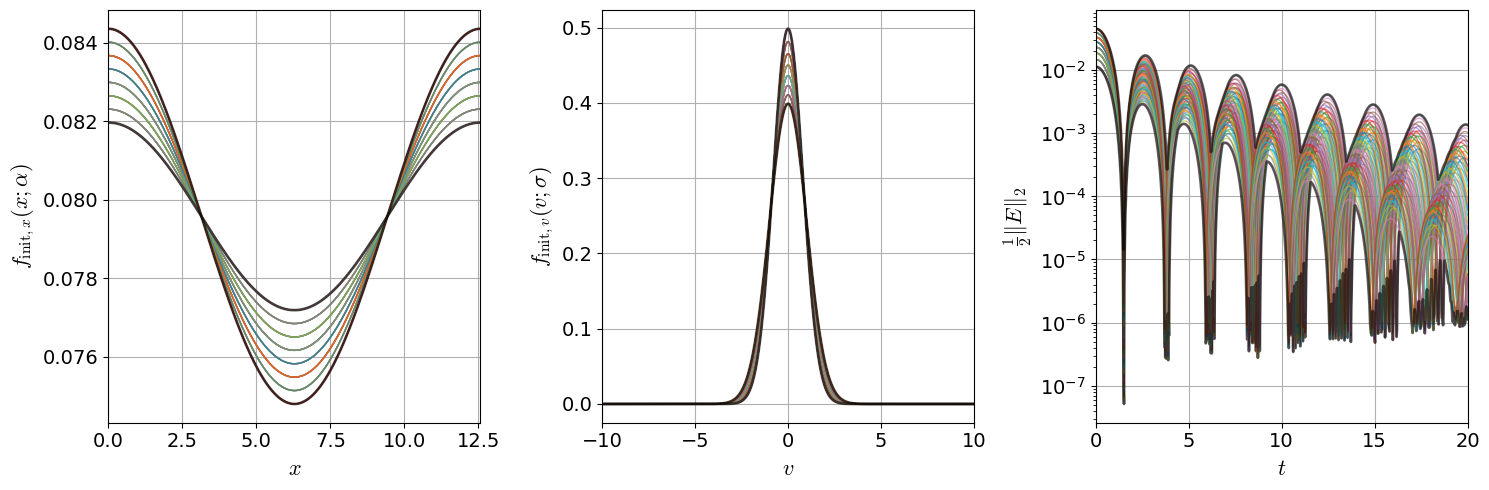}
    \caption{(Linear Landau damping) Initial distribution $f_{\init,x}(x;\alpha)$ (left), $f_{\init,v}(x;\sigma)$ (middle) and evolution of the electric energy $\ee\left(\bfX(t;\mu);\mu)\right)$ (right) for every $\mu \in \Gamma^{\train}$.}\label{fig: LLD_initial_dist_elec}
\end{figure}

We first consider the linear Landau damping test cases, with initial distributions \eqref{def:f_init_x}-\eqref{def:f_init_v_landaudamping} and $k=0.5$,
 $N = 10^5$ particles, $n_x = 48$ spatial discretization points. The final time equals $T = 20$ and the time step is set to $\Delta t = 2.5 \sci{-3}$. The parameter domain is taken equal to $\Gamma = [0.03, 0.06] \times [0.8, 1]$: the size of the perturbation is thus kept small. For the training dataset, we consider $P = 64$ parameters. The variation of the initial distribution and the electric energy damping as a function of $\mu$ is shown in \Cref{fig: LLD_initial_dist_elec}. Each color represents a different parameter in $\Gamma^{\train}$, and black lines are the envelopes of all the colored curves. 

\begin{figure}[htb!]
    \centering
    \includegraphics[width=\textwidth]{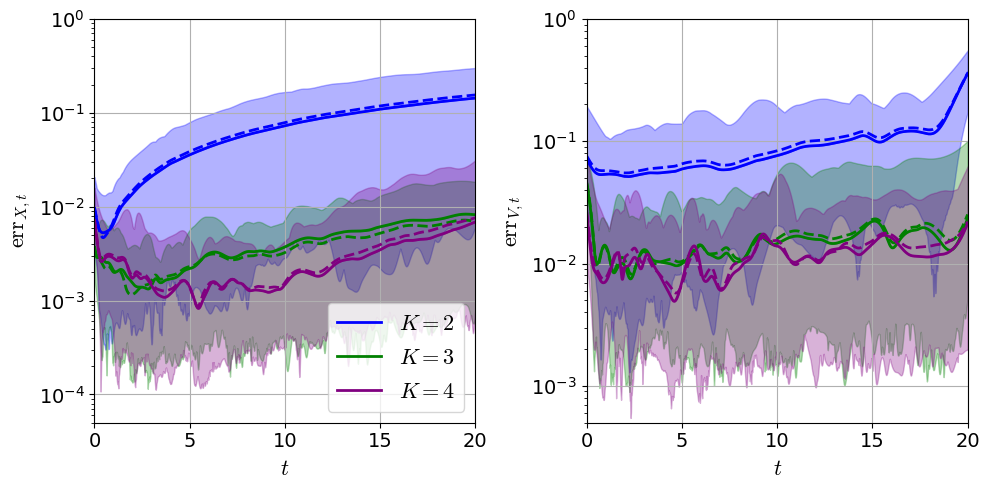}
    \caption{(Linear Landau damping) Mean error as a function of time $\error^{\opmean}_{X, t}$ (left,solid line) and $\error^{\opmean}_{V, t}$ (right, solid) for $\mu\in\Gamma^{\test}$. Each color stands for a value of $K\in\{2,3,4\}$, dashed lines are $\error^{\opmean}_{X, t}, \error^{\opmean}_{V, t}$ evaluated on the training set $\Gamma^{\train}$, the envelopes represent minimal and maximal errors $\error^{\opmin}_{X, t}, \error^{\opmax}_{X, t}$ (left) and  $\error^{\opmin}_{V, t}, \error^{\opmax}_{V, t}$.}\label{fig: LLD_error_Xt_Vt}
\end{figure}

The intermediate reduced variable size is set to $M=121$ and the complex SVD algorithm is used to build the first linear mapping. After completion of the training process with the architecture from \Cref{tab:hyperparameters pic hnn}, we evaluate our model on $\mu\in\Gamma^{\test}$. To begin with, we vary $K\in\{2,3,4\}$ and observe relative errors as a function of time in \Cref{fig: LLD_error_Xt_Vt}. A larger $K$ leads to smaller errors. For instance, with $K=2$, $\error^{\opmean}_{X, t}$ is around $2\sci{-2}$, while it is around $1\sci{-3}$ for $K=4$ . With this architecture, a reduced dimension $K=3$ is satisfactory. To obtain more precise results, we would have to modify the architecture presented in \Cref{tab:hyperparameters pic hnn} for a larger one. In the following, we set $K=3$.

\begin{figure}[htb!]
    \centering
    \includegraphics[width=\textwidth]{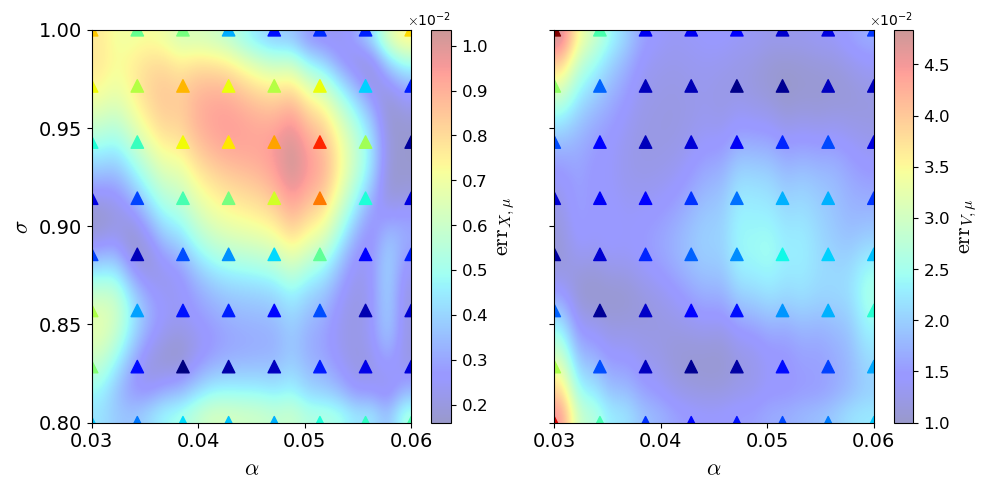}
    \caption{(Linear Landau damping) Errors as a function of the reduction parameters $\error_{X, \mu}$ (left) and $\error_{V, \mu}$ (right), triangular points represent the same error evaluated on the training set $\Gamma^{\train}$. }\label{fig: LLD_error_Xmu_Vmu_K3}
\end{figure}

In \Cref{fig: LLD_error_Xmu_Vmu_K3}, we then look at the relatives errors $\error_{X, \mu}, \error_{V, \mu}$ as a function of the parameters. The errors $\error_{X, \mu}, \error_{V, \mu}$ are of the order $6\sci{-3}$ and $3\sci{-2}$, respectively. In this specific case, we note that the maximal error is obtained inside the parameters domain for the positions and on the boundary for the velocities.

\begin{figure}[htb!]
    \centering
    \includegraphics[width=\textwidth]{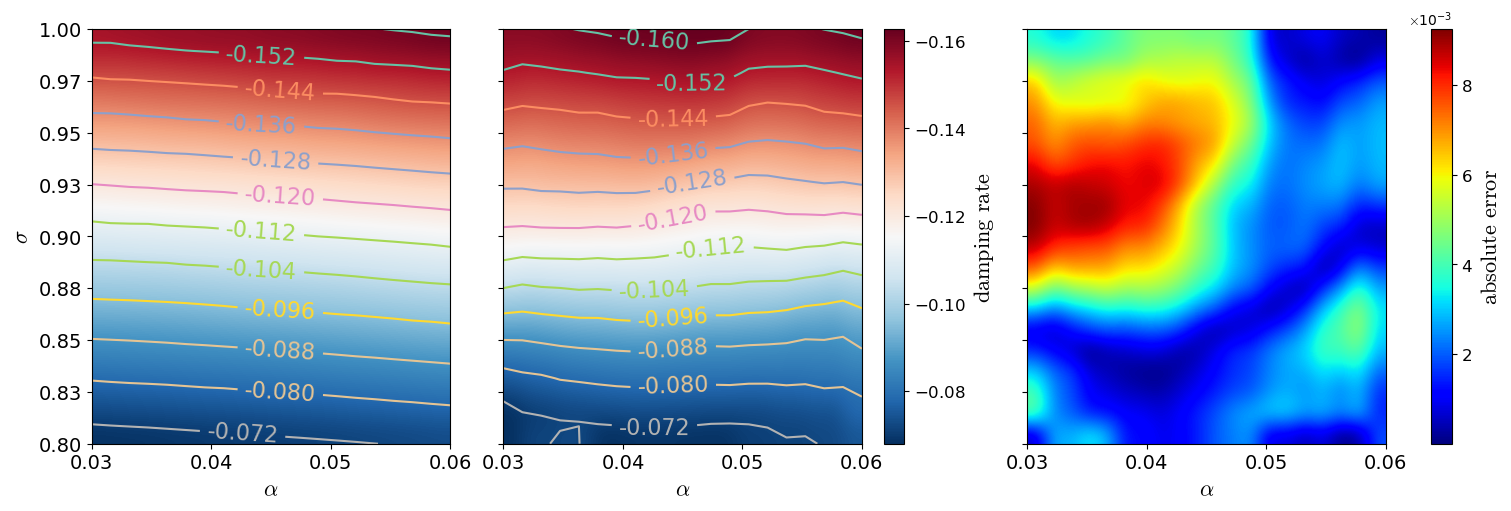}
    \caption{(Linear Landau damping) Electric energy $\ee\left(\bfX(t;\mu);\mu\right), \mu \in \Gamma^{\test}$ exponential damping rates of the FOM (left), the ROM (center) and absolute error (right).}\label{fig: LLD_error_damping_rate_K3}
\end{figure}

Next, we investigate the correctness of the damping rate. In theory, the electric energy $\ee\left(\bfX(t;\mu);\mu\right)$ decays exponentially in time with a constant damping rate, that depends on the standard deviation $\sigma$ of the Maxwellian initial distribution $f_{\init,v}$ and not on the amplitude $\alpha$ of the initial perturbation in space  $f_{\init, x}$ \cite{berge1969landau}. This property is captured by the reduced model, as observed in \Cref{fig: LLD_error_damping_rate_K3}. Thus, damping rates predictions are precise with an absolute error of about $5\sci{-3}$. In practice, the compression generates an error that causes the decay rate to fluctuate as a function of $\alpha$, but this fluctuation remains very small. Similarly, we can see that the reduced model slightly underestimates decay rates.

\begin{figure}[htb!]
    \begin{subfigure}{\textwidth}
        \centering
        \includegraphics[width=1.\textwidth]{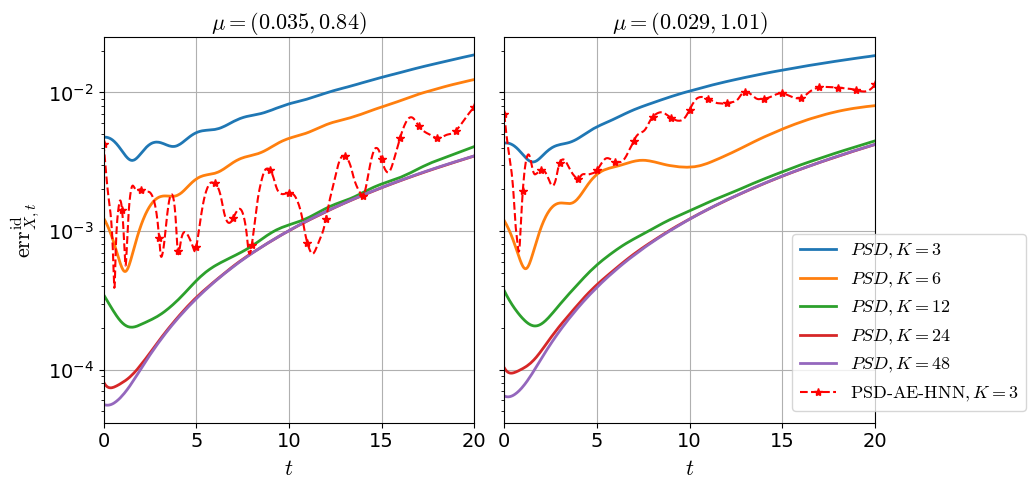}
        \caption{Errors $\error^{\text{id}}_{X, t}$}
        \label{fig: LLD_damping_psd_error_K3}
    \end{subfigure}
    \begin{subfigure}{\textwidth}
        \centering
        \includegraphics[width=1.\textwidth]{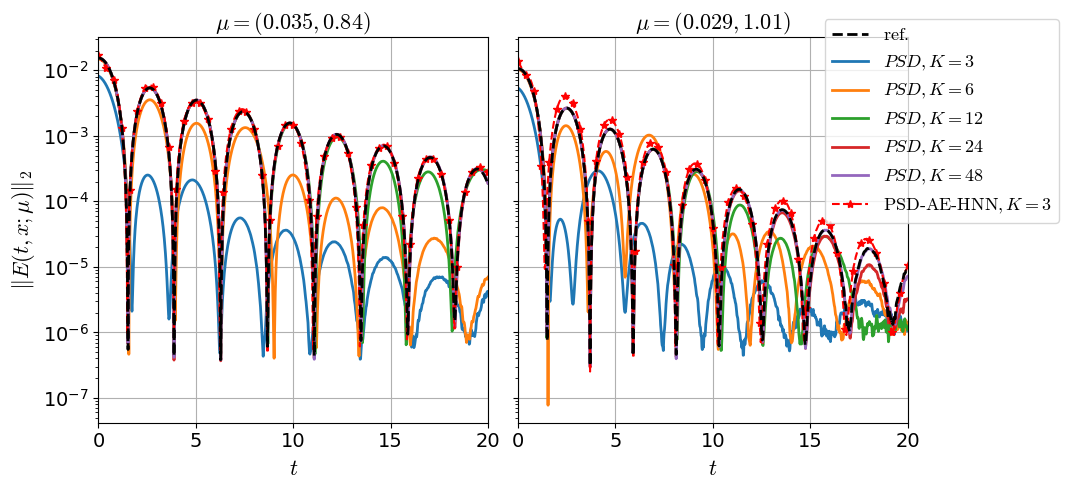}
        \caption{Electric energies $\ee\left(\bfX(t;\mu);\mu)\right)$}
        \label{fig: LLD_damping_psd_K3}
    \end{subfigure}
    \caption{(Linear Landau damping) Comparison of the PSD reduced model against our method with $K=3$, $\error^{\text{id}}_{X, t}=\left\| \bfX^t_\mu  - \hat\bfX^t_\mu \right\|_F / \left\| \bfX^t_\mu  \right\|_F$ for a given $\mu=(0.035, 0.84)$ (left) and $\mu=(0.029, 1.01)$ (right).}
    \label{fig: LLD_damping_psd_K3 and LLD_damping_psd_error_K3}
\end{figure}

We evaluate the performance of our method compared to the PSD-only approach in \Cref{fig: LLD_damping_psd_K3 and LLD_damping_psd_error_K3}. 
We test both methods with $K \in \{3, 6, 12, 24, 48\}$, and evaluate them at two parameter sets: $\mu = (0.035, 0.84) \in [0.03, 0.06] \times [0.8, 1]$ and $\mu = (0.029, 1.01) \notin [0.03, 0.06] \times [0.8, 1]$.
The damping rate, shown in \Cref{fig: LLD_damping_psd_K3}, indicates that $K \approx 30$ modes are required to match the performance of $K=3$ modes in the PSD-AE-HNN approach. However, as illustrated in \Cref{fig: LLD_damping_psd_error_K3}, the relative error in particle positions does not show significant differences. This highlights that the PSD method struggles to capture small-scale dynamics, that are crucial for preserving electric energy oscillations and damping, although it still performs well for large-scale dynamics.

    \subsection{Nonlinear Landau damping}\label{sec: non linear landau damping}

\begin{figure}[htb!]
    \centering
    \includegraphics[width=\textwidth]{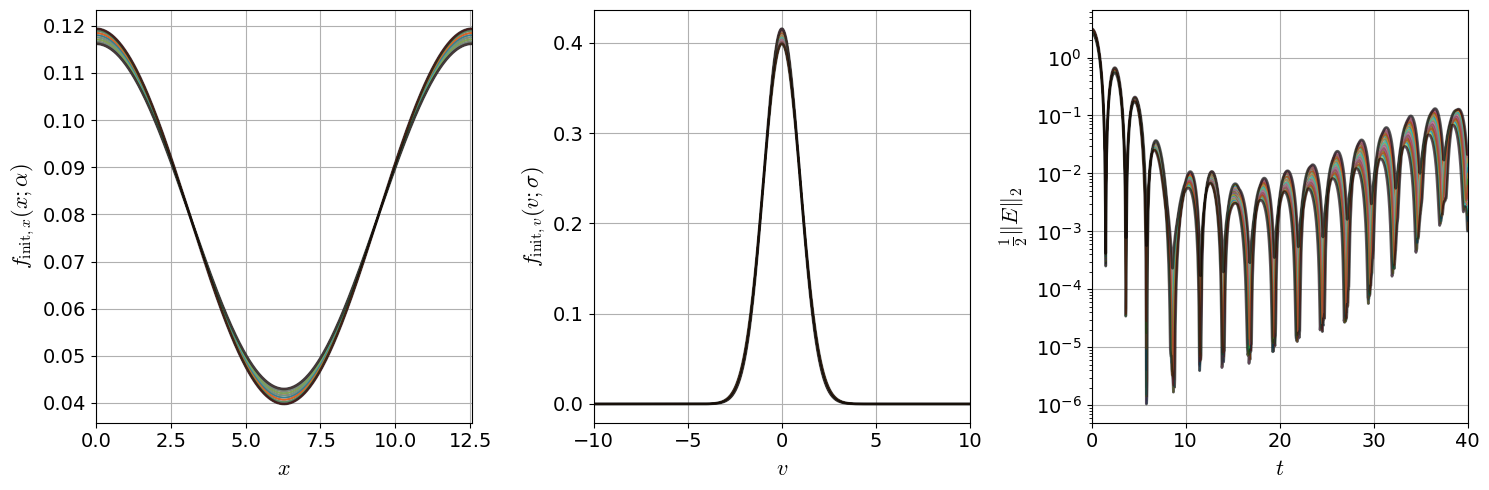}
    \caption{(Nonlinear Landau damping) Initial distribution $f_{\init,x}(x;\alpha)$ (left), $f_{\init,v}(x;\sigma)$ (middle) and evolution of the electric energy $\ee\left(\bfX(t;\mu);\mu\right)$ (right) for every $\mu \in \Gamma^{\train}$.}\label{fig: NLD_initial_dist_elec}
\end{figure}

In this test case, we keep the same initial distribution as in the previous section but consider a parametric domain, $\Gamma = [0.46,0.5] \times [0.96,1]$, with larger spatial perturbation amplitudes. We consider $\sci{5}$ particles, $n_x=64$ spatial discretization points. The final time and the time step are respectively set to  $T=40$ and $\Delta t=2.5\sci{-3}$. For the training dataset, we sample $(\alpha \; \sigma)^T \in \Gamma$ pairs over an $8 \times 8$ regular grid over $\Gamma$ forming $P=64$  pairs $\Gamma^{\train}$. In \Cref{fig: NLD_initial_dist_elec}, we plot the evolution of the initial distribution and the electric energy for all $\mu\in\Gamma^{\train}$. Each color represents a different value of $\mu$ and the envelopes are shown in black.

\begin{figure}[htb!]
    \centering
    \includegraphics[width=\textwidth]{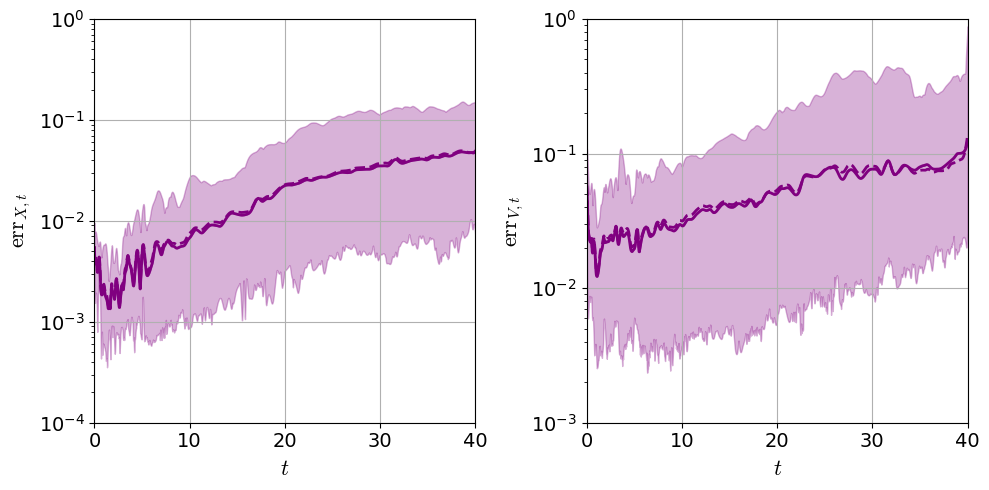}
    \caption{(Nonlinear Landau damping) Mean error as a function of time $\error^{\opmean}_{X, t}$ (left,solid line) and $\error^{\opmean}_{V, t}$ (right, solid) for $\mu\in\Gamma^{\test}$ for $K=4$. Dashed lines are $\error^{\opmean}_{X, t}, \error^{\opmean}_{V, t}$ evaluated on the training set $\Gamma^{\train}$, the envelopes represents minimal and maximal errors $\error^{\opmin}_{X, t}, \error^{\opmax}_{X, t}$ (left) and  $\error^{\opmin}_{V, t}, \error^{\opmax}_{V, t}$ (right).}\label{fig: NLD_error_Xt_Vt}
\end{figure}

Given the increased complexity compared with the linear Landau damping, the intermediate reduced dimension is set to $M=256$. The trained architecture is specified in \Cref{tab:hyperparameters pic hnn}. The reduced dimension is fixed equal to $K=4$. The relative errors as a function of time are depicted in \Cref{fig: NLD_error_Xt_Vt}. The errors in positions and velocities are both on the order of $1 \sci{-2}$.

\begin{figure}[htb!]
    \centering
    \includegraphics[width=\textwidth]{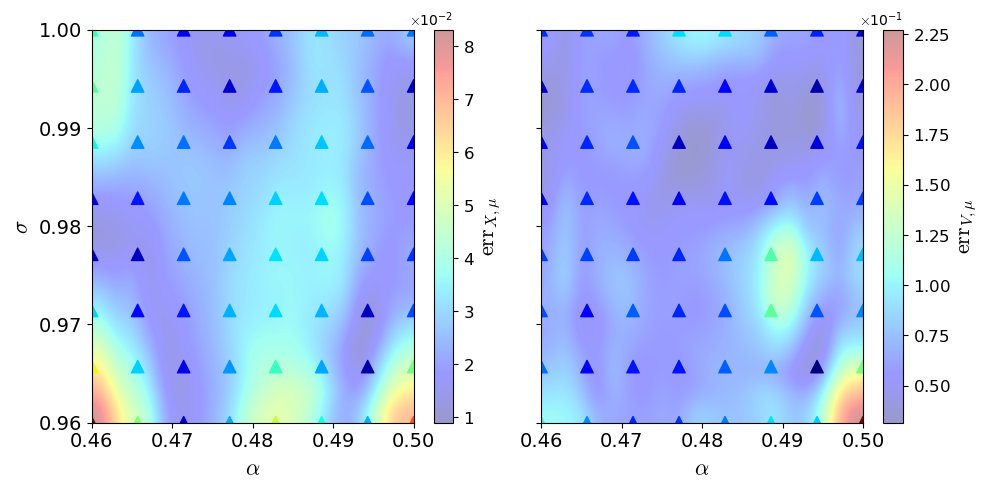}
    \caption{(Nonlinear Landau damping) Errors as a function of the reduction parameters $\error_{X, \mu}$ (left) and $\error_{V, \mu}$ (right), triangular points represent the same error evaluated on the training set $\Gamma^{\train}$. }\label{fig: NLD_error_Xmu_Vmu_K4}
\end{figure}

Subsequently, the relative errors as a function of $\mu$ are shown in \Cref{fig: NLD_error_Xmu_Vmu_K4}. The errors are around $2 \sci{-2}$ for the positions and $5 \sci{-2}$ for the velocities. As expected, errorr mainly concentrate on the boundary of the parameter domain.

\begin{figure}[htb!]
    \begin{subfigure}{\textwidth}
        \centering
        \includegraphics[width=1.\textwidth]{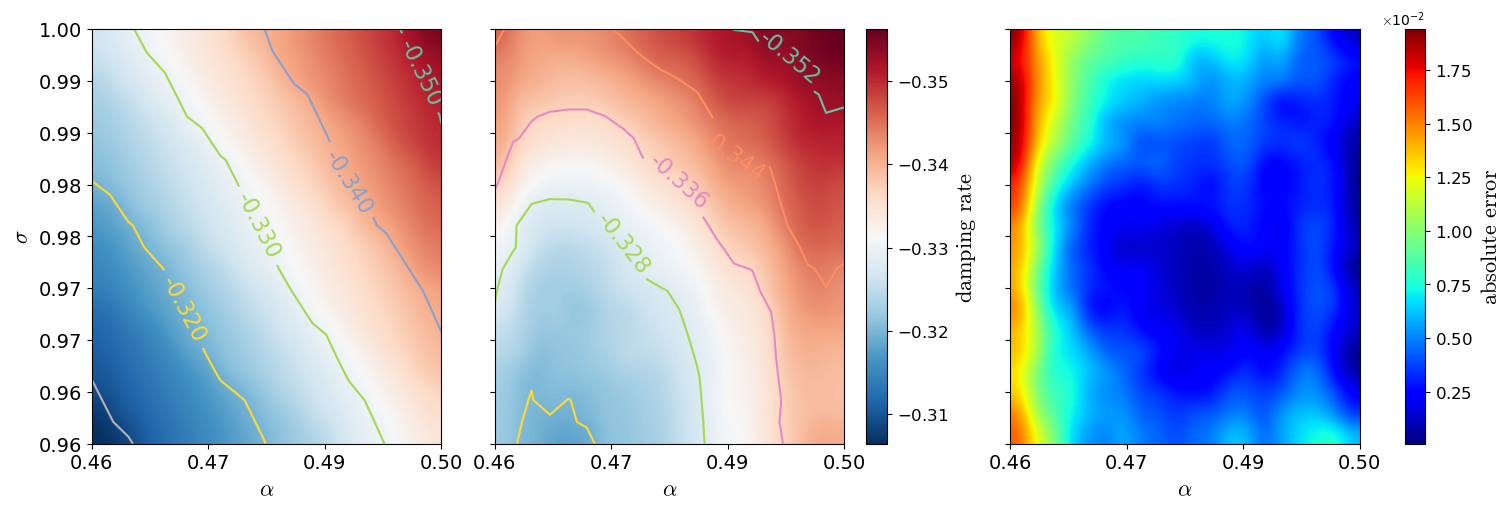}
        \caption{Damping rate.}
        \label{fig: NLD_error_damping_rate_K4}
    \end{subfigure}
    \begin{subfigure}{\textwidth}
        \centering
        \includegraphics[width=1.\textwidth]{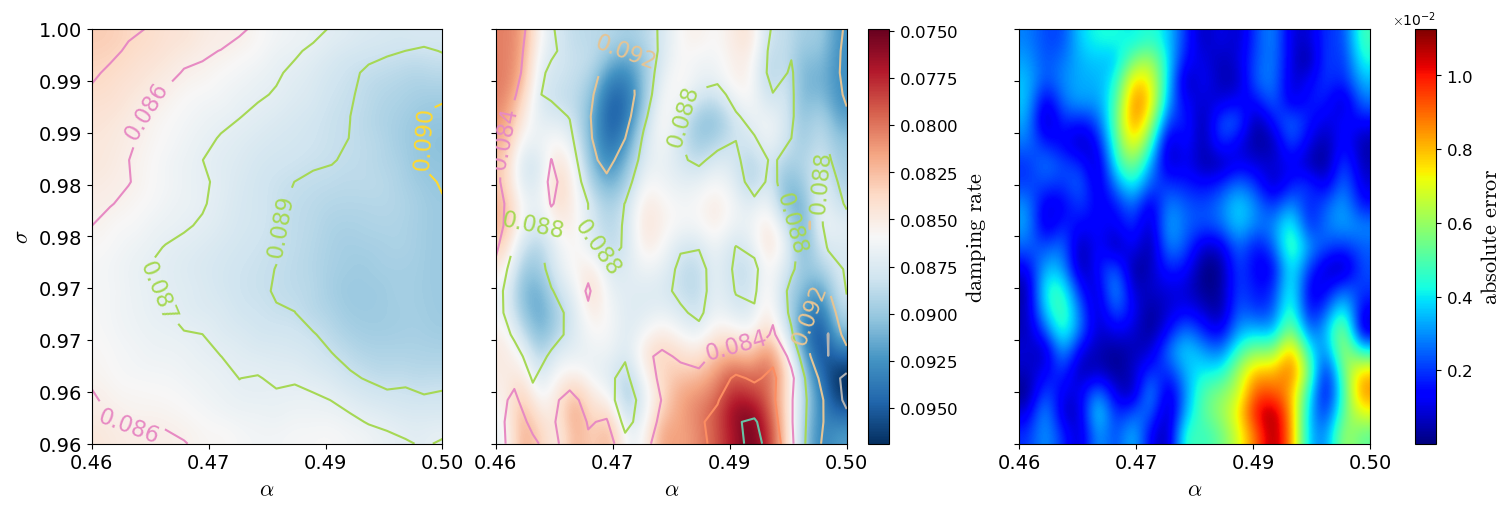}
        \caption{Growth rate.}
        \label{fig: NLD_error_growth_rate_K4}
    \end{subfigure}
    \caption{(Nonlinear Landau damping) Electric energy $\ee\left(\bfX(t;\mu);\mu\right), \mu \in \Gamma^{\test}$ exponential damping rates (top) and growth rates (bottom) of the FOM (left), the ROM (center) and absolute error (right).}
    \label{fig: NLD_error_damping_rate_K4 and NLD_error_growth_rate_K4}
\end{figure}

Then, in \Cref{fig: NLD_error_damping_rate_K4 and NLD_error_growth_rate_K4}, we compare the exponential damping and growth rates of the electric energy in ROM with those in FOM.. As observed in \Cref{fig: NLD_initial_dist_elec}, we expect a constant damping rate when $t<10$ then a constant growth rate when $t>20$. \Cref{fig: NLD_error_damping_rate_K4} shows that the mean error on the damping rate is about $7\sci{-3}$ and the error is maximal for the smallest values of $\alpha$. \Cref{fig: NLD_error_growth_rate_K4} shows that the mean error on the growth rate is about $3\sci{-3}$. On the other hand, unlike the linear case, the dependency of the rates to the two parameters is less well captured.

\begin{figure}[htb!]
    \centering
    \includegraphics[width=\textwidth]{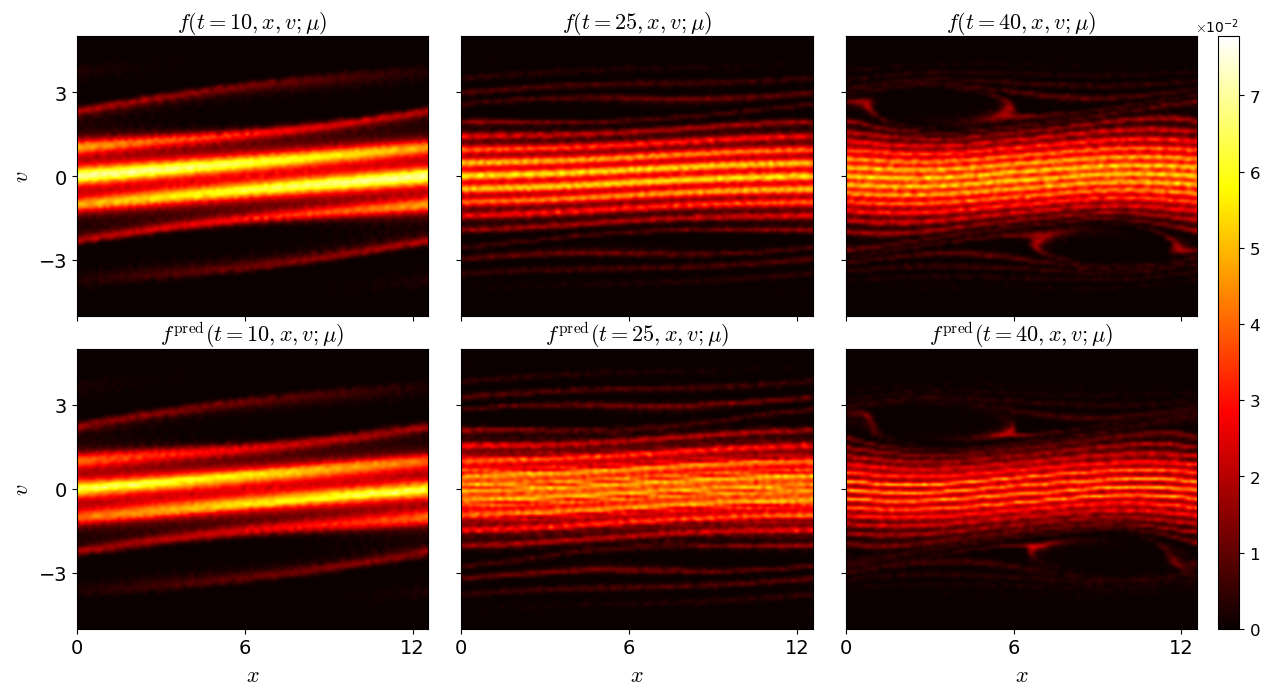}
    \caption{(Nonlinear Landau damping) Solution $f(t,x,v;\mu)$ (top) and $f^{\text{pred}}(t,x,v;\mu)$ (bottom) for $t\in\{10,25,40\}$ and $\mu=(0.465, 0.986)$. }\label{fig: NLD_distribution_plot}
\end{figure}

Next, we compare the evolution of the distribution $f(t,x,v;\mu)$ from the reference model with the predictions of its reduced model in \Cref{fig: NLD_distribution_plot}. While small differences are observed, the overall dynamics are well captured.

\begin{figure}[htb!]
    \begin{subfigure}{\textwidth}
        \centering
        \includegraphics[width=1.\textwidth]{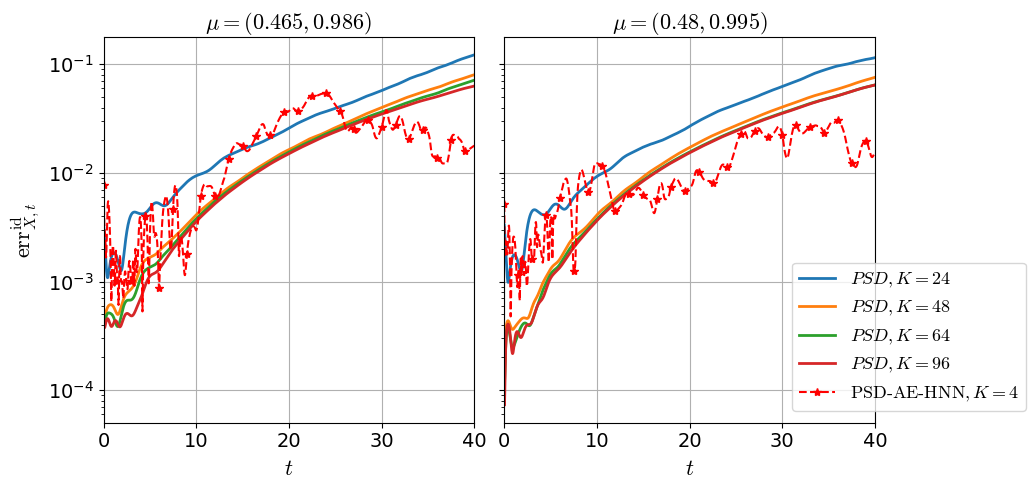}
        \caption{Errors $\error^{\text{id}}_{X, t}$}
        \label{fig: NLD_damping_psd_error_K4}
    \end{subfigure}
    \begin{subfigure}{\textwidth}
        \centering
        \includegraphics[width=1.\textwidth]{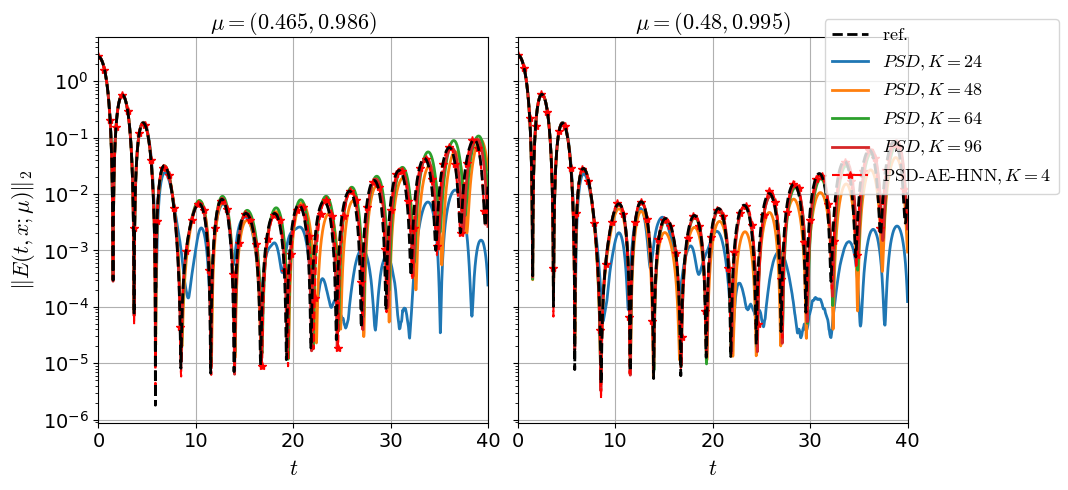}
        \caption{Electric energies $\ee\left(\bfX(t;\mu);\mu\right)$}
        \label{fig: NLD_damping_psd_K4}
    \end{subfigure}
    \caption{(Nonlinear Landau damping) Comparison of the PSD reduced model against our method with $K=4$, $\error^{\text{id}}_{X, t}=\left\| \bfX^t_\mu  - \hat\bfX^t_\mu \right\|_F / \left\| \bfX^t_\mu  \right\|_F$ for a given $\mu=(0.465, 0.986)$ (left) and $\mu=(0.48, 0.995)$ (right).}
    \label{fig: NLD_damping_psd_K4 and NLD_damping_psd_error_K4}
\end{figure}

In \Cref{fig: NLD_damping_psd_K4 and NLD_damping_psd_error_K4}, we then compare the precision of the method with the PSD-only approach for two test parameters $\mu=(0.465, 0.986)$ and $\mu=(0.48, 0.995)$. To match the performance of our method, about $K=100$ PSD modes are required, both for the particle distribution errors (\Cref{fig: NLD_damping_psd_error_K4}) and for the electric energy (\Cref{fig: NLD_damping_psd_K4}). In this test case, our approach is particularly effective for the positions and velocities associated with the last oscillations of the simulation.

\begin{figure}[htb!]
    \centering
    \includegraphics[width=\textwidth]{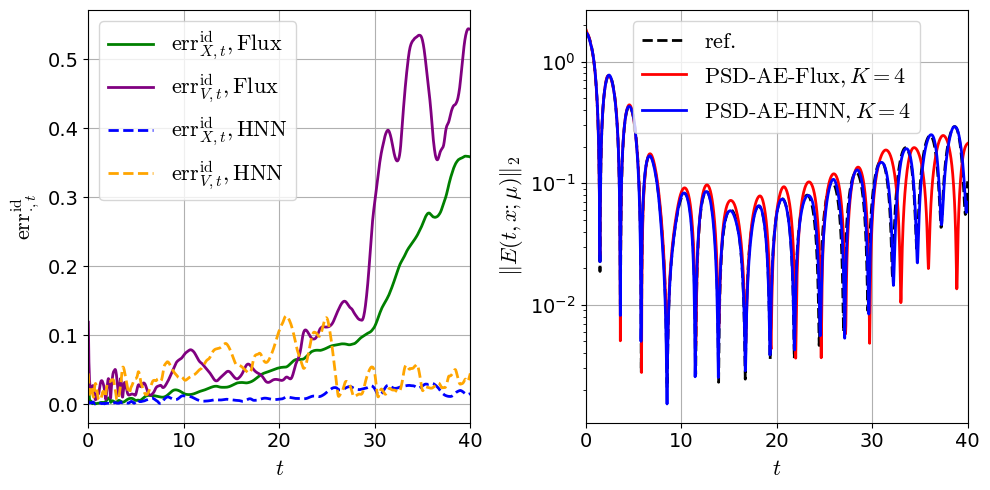}
    \caption{(Nonlinear Landau damping) PSD-AE-Flux prediction for a single test parameter $\mu$ compared to the PSD-AE-HNN method. Errors as a function of time $\error^{\text{id}}_{X, t}, \error^{\text{id}}_{V, t}$ (left) and predicted electric energy $\ee\left(\bfX(t;\mu);\mu)\right)$.}\label{fig: NLD_baseline_output_test1}
\end{figure}

Finally, we assess the importance of using a HNN in our PSD-AE-HNN method, compared to a flux-approximating neural approach. For this, we replace the HNN in \Cref{sec: A Hamiltonian reduction with Proper Symplectic Decompostion prereduction} with a standard neural network $\bar{\mathcal{F}}_{\theta_f}$, which approximates the flux of the reduced model:
\begin{align*}\label{eq:baseline rom}
    \begin{cases}
    \displaystyle \frac{d}{dt} \bfu(t;\mu) = \bar{\mathcal{F}}_{\theta_f}\left(\bfu(t;\mu)\right) , &\text{ in } [0,T],\\
    \bfu(0;\mu) = \mathcal{E}\left(\bfU_{\init}(\mu)\right).
   \end{cases}
\end{align*} 
We also replace, in the prediction operator $\mathcal{P}_s$, the symplectic Störmer-Verlet by the standard Runge-Kutta 2 scheme. The stability loss weight is set to $\omega_{\stab} = 0$ as it cannot be evaluated and the remaining parameters are unchanged. The trained model is tested on $\mu = (0.48, 0.995) \in \Gamma$, and the results are shown in \Cref{fig: NLD_baseline_output_test1}. The errors increase strongly over time, and the electric energy growth is not accurately captured—unlike our PSD-AE-HNN method, which performs accurately in comparison. Unlike an approach based solely on PSD, where the reduced model is Hamiltonian and the encoder and decoder are symplectic, here only the reduced model is Hamiltonian. This last case shows that simply preserving the structure in the reduced space is sufficient to improve long-term stability.

    \subsection{Two-stream instability}\label{sec: two-stream instability}

\begin{figure}[htb!]
    \centering
    \includegraphics[width=\textwidth]{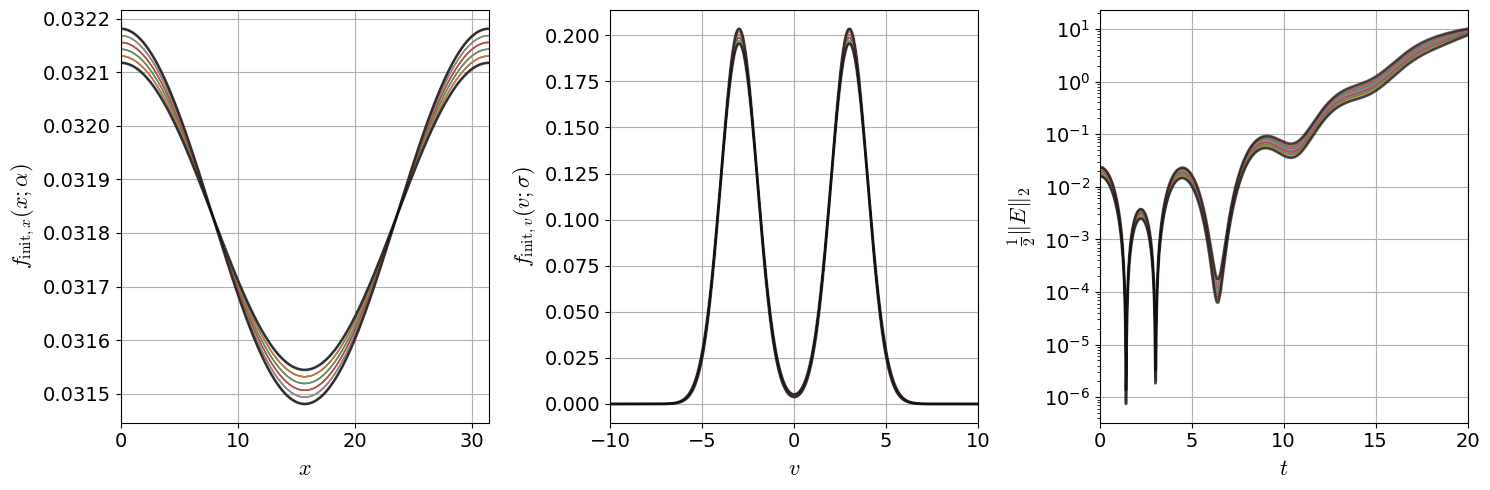}
    \caption{(Two stream instability) Initial distribution $f_{\init,x}(x;\alpha)$ (left), $f_{\init,v}(x;\sigma)$ (middle) and evolution of the electric energy $\ee\left(\bfX(t;\mu);\mu)\right)$ (right) for every $\mu \in \Gamma^{\train}$.}\label{fig: TSI_initial_dist_elec}
\end{figure}

In this third test case, we consider the initial distributions \eqref{def:f_init_x}-\eqref{def:f_init_v_twostream}, with the sum of two Gaussians in velocity, and $k=0.2$. The number of particles equals  $N = 1.5 \sci{5}$ and there are $n_x=64$ spatial discretization points. The final time equals $T=20$ and the time step $\Delta t = 2.5 \sci{-3}$. The parameter domain is taken equal to $\Gamma = [0.009, 0.011] \times [0.98, 1.02]$ and the training set is composed $P=36$ distinct pairs. This training set is more scattered compared to the other test cases. For comparison purposes, the authors in \cite{hesthaven2024adaptive} use $300$ snapshots for the same test case. The variation in the initial distribution and electric energy is shown in \Cref{fig: TSI_initial_dist_elec}.

\begin{figure}[htb!]
    \centering
    \includegraphics[width=\textwidth]{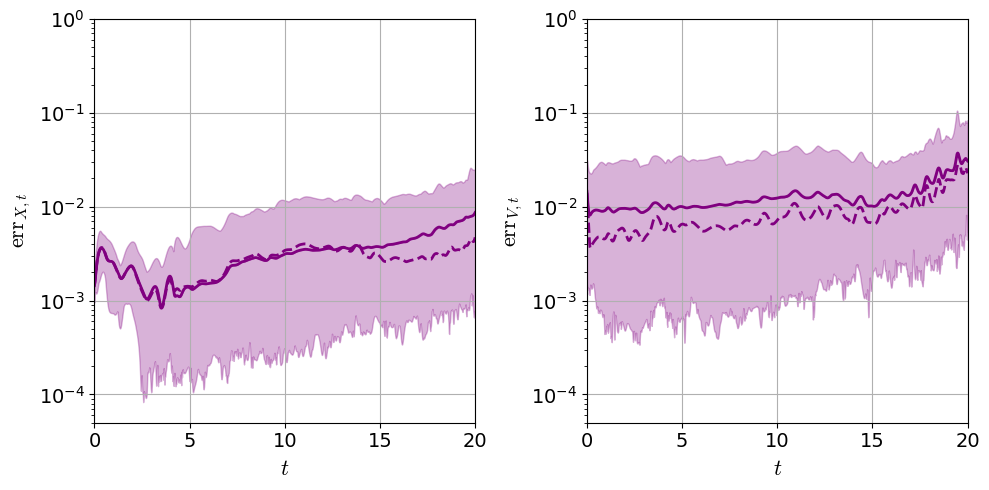}
    \caption{(Two stream instability) Mean error as a function of time $\error^{\opmean}_{X, t}$ (left,solid line) and $\error^{\opmean}_{V, t}$ (right, solid) for $\mu\in\Gamma^{\test}$. Dashed lines are $\error^{\opmean}_{X, t}, \error^{\opmean}_{V, t}$ evaluated on the training set $\Gamma^{\train}$, the envelopes represents minimal and maximal errors $\error^{\opmin}_{X, t}, \error^{\opmax}_{X, t}$ (left) and  $\error^{\opmin}_{V, t}, \error^{\opmax}_{V, t}$.}\label{fig: TSI_error_Xt_Vt}
\end{figure}

The intermediate reduced dimension and the final reduced dimensions are set to $M=121$ and $K=4$. The AE-HNN networks is trained with the architecture presented in \Cref{tab:hyperparameters pic hnn}. We inspect relatives errors as a function of time in \Cref{fig: TSI_error_Xt_Vt}: $\error^{\opmean}_{X, t}$ is about $2\sci{-3}$ and $\error^{\opmean}_{V, t}$ is on the order of $1\sci{-2}$.

\begin{figure}[htb!]
    \centering
    \includegraphics[width=\textwidth]{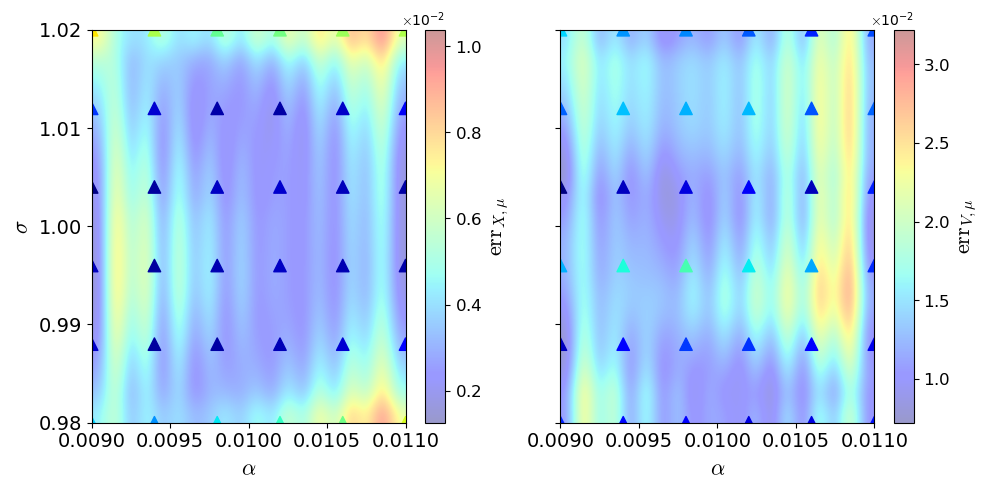}
    \caption{(Two stream instability) Errors as a function of the reduction parameters $\error_{X, \mu}$ (left) and $\error_{V, \mu}$ (right), triangular points represent the same error evaluated on the training set $\Gamma^{\train}$. }\label{fig: TSI_error_Xmu_Vmu_K4}
\end{figure}

Next, we observe the relative errors as a function of $\mu$ in \Cref{fig: TSI_error_Xmu_Vmu_K4}. The error in positions, $\error_{X, \mu}$, is approximately $5\sci{-3}$, and the error in velocities, $\error_{V, \mu}$, is around $2\sci{-2}$. We notice only a slight increase in error  as $\alpha$ increases, attributed to the sparsity of the training set. Overall, errors remain low, and the reduced dynamics are learned effectively.

\begin{figure}[htb!]
    \centering
    \includegraphics[width=\textwidth]{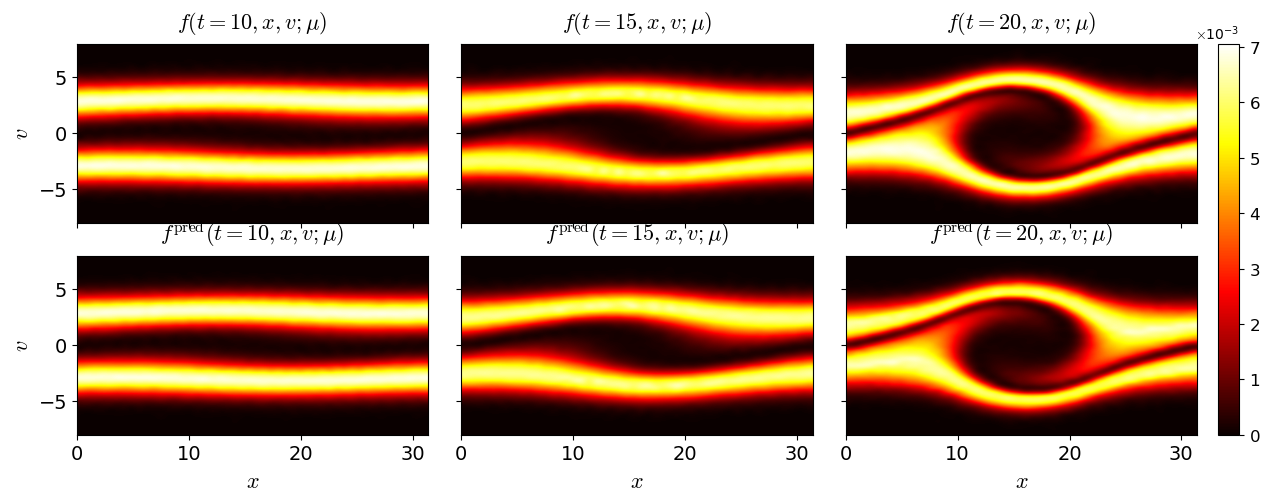}
    \caption{(Two stream instability) $f(t,x,v;\mu)$ (top) and $f^{\text{pred}}(t,x,v;\mu)$ (bottom) for $t\in\{10,15,20\}$ and $\mu=(0.0095, 0.99)$. }\label{fig: TSI_distribution_plot}
\end{figure}

Then, we plot the evolution of the distribution $f(t,x,v;\mu)$ of the FOM at different times in comparison with the ROM predicted distribution $f^{\text{pred}}(t,x,v;\mu)$ in \Cref{fig: TSI_distribution_plot}. The dynamics are correctly captured from the initial stream shearing to the development of a central vortex. 

\begin{figure}[htb!]
    \begin{subfigure}{\textwidth}
        \centering
        \includegraphics[width=1.\textwidth]{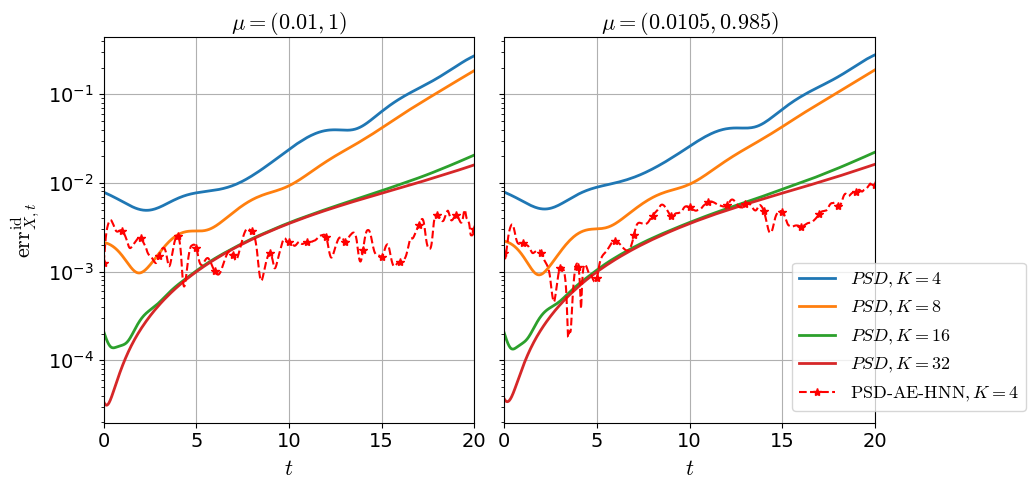}
        \caption{Errors $\error^{\text{id}}_{X, t}$}
        \label{fig: TSI_damping_psd_error_K4}
    \end{subfigure}
    \begin{subfigure}{\textwidth}
        \centering
        \includegraphics[width=1.\textwidth]{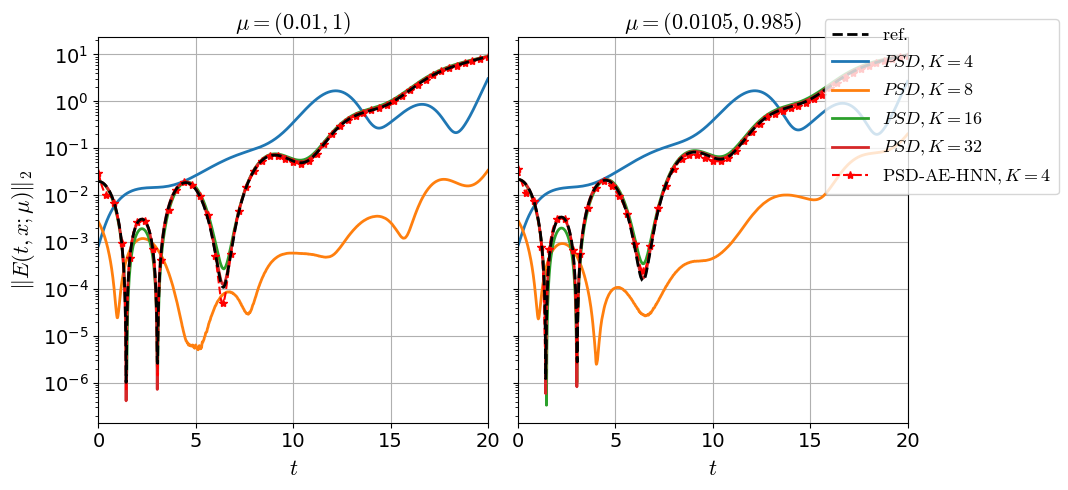}
        \caption{Electric energies $\ee\left(\bfX(t;\mu);\mu)\right)$}
        \label{fig: TSI_damping_psd_K4}
    \end{subfigure}
    \caption{(Two stream instability) Comparison of the PSD reduced model against our method with $K=4$, $\error^{\text{id}}_{X, t}=\left\| \bfX^t_\mu  - \hat\bfX^t_\mu \right\|_F / \left\| \bfX^t_\mu  \right\|_F$ for a given $\mu=(0.01, 1)$ (left) and $\mu=(0.0105, 0.985)$ (right).}
    \label{fig: TSI_damping_psd_K4 and TSI_damping_psd_error_K4}
\end{figure}

Eventually, we compare the method with $K=4$ to the PSD with $K\in\{4,8,16,32\}$ in \Cref{fig: TSI_damping_psd_K4 and TSI_damping_psd_error_K4} for $\mu=(0.01, 1)$ (left) and $(0.0105, 0.985)$ (right). In \Cref{fig: TSI_damping_psd_K4}, we observe the electric energy $\ee\left(\bfX(t;\mu);\mu\right)$, where its dynamics are well replicated with our method. We would need $K=30$ modes with the PSD to obtain comparable results. In \Cref{fig: TSI_damping_psd_error_K4}, we study the relative errors $\error^{\text{id}}_{X, t}$ and conclude that our method with a dimension of $K=4$ achieves comparable results in terms of precision with the PSD and $K=30$ modes.

\subsection{Computation gain}
    
In this section, we evaluate the computation time for each model and test case from \Cref{sec: linear landau damping,sec: non linear landau damping,sec: two-stream instability}. To compare the performance of the reduced model, it is essential to identify a full model with equivalent accuracy. To achieve this, the number of particles will be varied. This approach will help us to assess that our method offers superior performance compared to a simple reduction in the number of particles. Therefore, we will not discuss the execution time of the PSD-only reduced model described in \Cref{eq: psd reduced model}, as it is expected to result in a longer execution time without any hyper-reduction techniques.

In this study, we focus on a specific quantity of interest for each test case and evaluate the computation time for a single parameter, $\mu$, across varying particle numbers, $N$. It is important to note that PIC simulations tend to exhibit noise when the particle count is low, leading to a non-monotonic relationship between the error and the particle number, $N$, with respect to the quantity of interest.

The code is implemented in Python, with the majority of operations utilizing the NumPy library. However, neural networks and training processes are implemented using the TensorFlow framework. A single AMD Ryzen 9 5900X CPU is employed for computation. 
It is important to note that this runtime test does not fully reflect the strengths of our method for two main reasons: (i) it is run on a CPU rather than a GPU, which significantly limits the efficiency of neural network evaluations; and (ii) it considers only a single initial condition, which prevents us from showcasing the neural network’s ability to handle multiple inputs efficiently through vectorization.

\begin{figure}[htb!]
    \centering
    \begin{subfigure}[b]{1.\textwidth}
        \centering
        \includegraphics[width=0.75\textwidth]{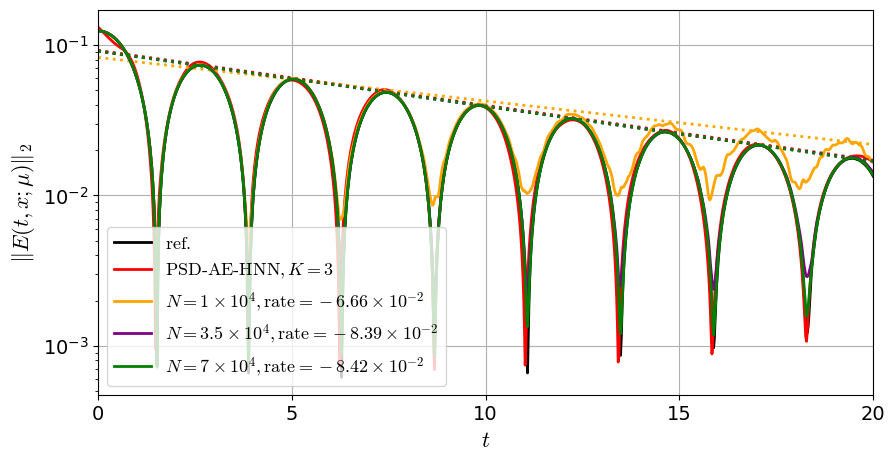}
    \caption{(Linear Landau damping) $\mu=(0.035, 0.84)$. Dashed lines represent the linear damping.}\label{fig: LLD_quantity_of_interest_N}
    \end{subfigure}

    \begin{subfigure}[b]{1.\textwidth}
        \centering
        \includegraphics[width=0.75\textwidth]{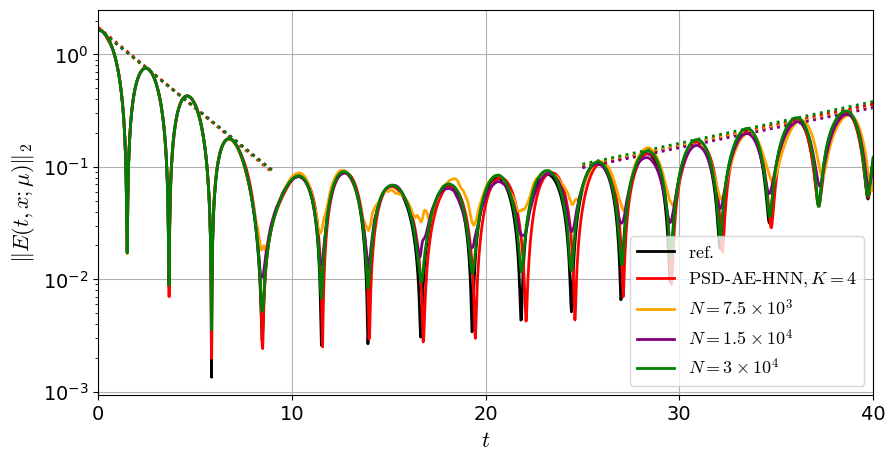}
    \caption{(Nonlinear Landau damping) $\mu=(0.465, 0.986)$. Dashed lines represent the linear damping and growth.}\label{fig: NLD_quantity_of_interest_N}
    \end{subfigure}

    \begin{subfigure}[b]{1.\textwidth}
        \centering
        \includegraphics[width=0.75\textwidth]{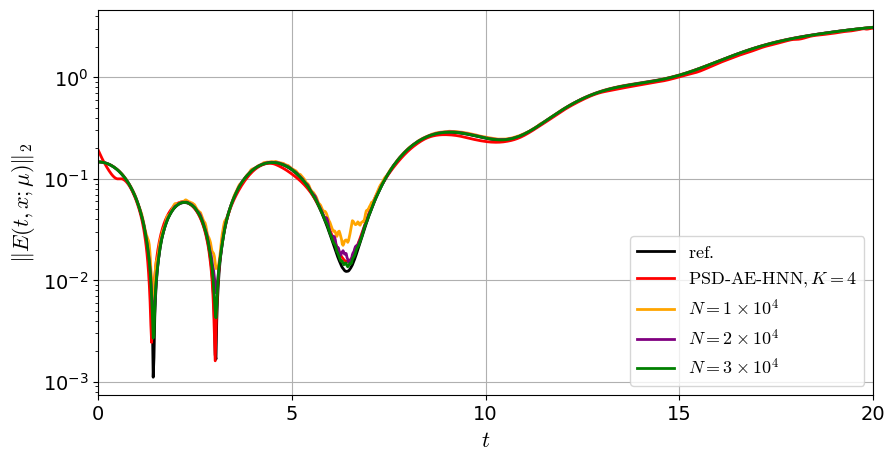}
    \caption{(Two-stream instability) $\mu=(0.0105, 0.985)$.}\label{fig: TSI_quantity_of_interest_N}
    \end{subfigure}
    \caption{Electric energy as a function of time for various $N$ and the PSD-AE-HNN reduced model. Solid lines represent the energies.}
    \label{fig: quantity_of_interest_N}
\end{figure}

For the linear Landau damping from \Cref{sec: linear landau damping}, we set $\mu=(0.035, 0.84)$. We focus on the estimation of the damping rate to test our method. As shown in \Cref{fig: LLD_quantity_of_interest_N}, we estimate a PIC simulation  with $N=7\sci4$ particles to be as precise as our PSD-AE-HNN method. From \Cref{tab: LLD_speed}, the PSD-AE-HNN is $4.63$ times faster than the prior.

Then, we focus on the nonlinear Landau damping test case with $\mu=(0.465, 0.986)$. We consider the damping and growth rates as the quantities of interest. As shown in \Cref{fig: NLD_quantity_of_interest_N}, the equivalent PIC simulation requires  $N=3\sci{4}$ particles and the speedup is about $1.95$.

Finally, we are interested in the evolution of the electric energy for the two-stream instability. We set $\mu=(0.0105, 0.985)$ and observe that we need about $N=3\sci{4}$ particles for a comparable precision to our method in \Cref{fig: TSI_quantity_of_interest_N} giving an acceleration of $2.10$.

\begin{figure}[htb!]
    \centering
    \begin{subfigure}[b]{1.\textwidth}
        \centering
        \begin{tabular}{cccccc}
        \toprule
&\multicolumn{4}{c}{PIC} & PSD-AE-HNN \\
N & $1\sci{5}$ & $7\sci{4}$ & $3.5\sci{4}$ & $1\sci{4}$ &\\ \midrule
damping rate ($\sci{-2}$)& $-8.44$ & $-8.42$ & $-8.39$ & $-6.66$ & $-8.41$ \\ \midrule
time (s) & 25.05 & 11.40 & 6.13 & 2.00 & 2.46 \\ \midrule
speedup & 0.46 & 1 & 1.86 & 5.70 & 4.63 \\ \bottomrule
\end{tabular}\caption{Linear Landau damping}
        \label{tab: LLD_speed}
    \end{subfigure}
    \vspace{0.7cm} 

    \begin{subfigure}[b]{1.\textwidth}
        \centering
\begin{tabular}{cccccc}
\toprule
&\multicolumn{4}{c}{PIC} & PSD-AE-HNN \\
N & $1\sci{5}$ & $3\sci{4}$ & $1.5\sci{4}$ & $7.5\sci{3}$ & \\ \midrule
damping rate ($\sci{-1}$) & -3.23 & -3.23 & -3.23 & -3.26 & -3.31 \\ \midrule
growth rate ($\sci{-2}$) & 8.55 & 8.55 & 8.22 & 7.89 & 8.60 \\ \midrule
time (s) & 53.45 & 11.13 & 6.03 & 3.37 & 5.71 \\ \midrule
speedup & 0.21 & 1 & 1.85 & 3.30 & 1.95 \\ \bottomrule
\end{tabular}\caption{Nonlinear Landau damping}
        \label{tab: NLD_speed}
    \end{subfigure}
    \vspace{0.7cm} 

    \begin{subfigure}[b]{1.\textwidth}
        \centering
\begin{tabular}{cccccc}
\toprule
&\multicolumn{4}{c}{PIC} & PSD-AE-HNN \\
N & $1.5\sci{5}$ & $3\sci{4}$ & $2\sci{4}$ & $1\sci{4}$ & \\ \midrule
time (s) & 59.30 & 5.34 & 3.68 & 2.08 & 2.54 \\ \midrule
speedup & 0.09 & 1 & 1.45 & 2.57 & 2.10 \\ \bottomrule
\end{tabular}\caption{Two-stream instability}
        \label{tab: TSI_speed}
    \end{subfigure}
    \caption{Computation time and numerical acceleration for each test case.}
    \label{fig: speed_per_test_case}
\end{figure}

From a theoretical point of view and discarding the projection cost, one integration step requires $O(N + n_x^2)$ operations. From \cite{cote2025hamiltonian}, our reduced model cost is about $O(\sum_{k=1}^L n^{(k - 1)}n^{(k)})$ where $n^{(k)}$ is $n$-th layer width i.e. number of units of the HNN. In addition, as $n^{(k)}$ is of the order of $K^2$, the cost is $O(K^4)$ which makes it competitive as it does not depend on $N$ nor $n_x$.


\section{Conclusion}

We have introduced a new Hamiltonian reduction method to reduce the number of particles in a particle-based discretization of the Vlasov-Poisson equation. This method uses a two-step mapping that combines Proper Symplectic Decomposition (PSD) for linear reduction with an autoencoder (AE) for additional nonlinear compression. PSD significantly reduces the dimensionality of the problem, while AE further compresses the data in a nonlinear manner. The reduced dynamics are then captured by a Hamiltonian neural network (HNN), ensuring the preservation of the Hamiltonian structure. The PSD-AE-HNN approach shows strong performance in linear and nonlinear test cases compared to the PSD method and offers good computational efficiency. Overall, the method is data-driven, non-intrusive, and highly adaptable.

One issue that arises is how to improve the quality of the approximation. In projection-based methods such as PSD, increasing the reduced dimension $K$ leads to an increase in accuracy. There is no systematic process of this type for the proposed method. Instead, improvements can come from tuning the hyperparameters of the neural networks (e.g., the number of layers, size of the layers, activation functions) and refining the learning strategy.

The results should be extended to cover PIC simulations in two or three spatial and velocity dimensions. A substantial increase in network size should not be necessary for the AE-HNN component, as convolutional layers can be used. In contrast, the PSD may require further refinement. Additionally, future extensions could also involve integrating time-adaptive model reduction techniques as proposed in \cite{hesthaven2022reduced}.

\paragraph{Acknowledgements.} This research was funded in part by l’Agence Nationale de la Recherche (ANR), project
ANR-21-CE46-0014 (Milk).


\appendix

\section{A short overview of the PSD}\label{annex : PSD details}

Here we make a more detailed presentation of the Proper Symplectic Decomposition (PSD) \cite{peng2015symplectic}. For this method, the PSD decoder is a linear symplectic operator, as presented in Sections ~\ref{sec: PSD and AE-HNN coupled reduced order model} and \ref{sec: psd details}. It writes
\[ 
    \widetilde{\bfU} = A^+ \bfU, 
\]
with $A \in  \text{Sp}_{2K,2N}(\R)$ a symplectic matrix. This matrix is determined by minimizing the reconstruction error of the projection:
\[
    \underset{A \in \text{Sp}_{2N,2K}(\R)}{\operatorname{min}}\  \| U - A A^+ U \|_F,
\]
where $U$ refers to the snapshot matrix defined in Eq.~\ref{eq: U snaphots matrix}. For a direct computation of $A$, we search for a solution in the following subset of symplectic matrix
\[
    \text{Sp}_{2N,2K}(\R) \,\cap\, \left\{ \begin{pmatrix}
        \Phi & - \Psi \\ \Psi & \Phi
    \end{pmatrix}, \Phi, \Psi \in \mathcal{M}_{N,K}(\R) \right\},
\]
using the Complex SVD method. 
The snapshot matrix is transformed into
\begin{equation*}
    U' = \left\{
  \bfX^0_{\mu_1}, \dots, \bfX^{n_T}_{\mu_P}
    \right\} + 
    i \, \left\{
  \bfV^0_{\mu_1}, \dots, \bfV^{n_T}_{\mu_P}
    \right\} \in \mathcal{M}_{N, (n_T + 1)P}(\R)(\mathbb{C}).
\end{equation*}
The minimization problem reads
\begin{equation}\label{eq: min problem PSD 2}
    \underset{B \in \operatorname{U}_{N,M}(\mathbb{C}), B^{\ast} B=I_M}{\operatorname{min}} \left\| U' - BB^\ast U' \right\|_F,
\end{equation}
where $\operatorname{U}_{N,M}(\mathbb{C})$ denotes the set of unitary matrices and, for any matrix $B$, $B^\ast$ is its conjugate-transpose. As a consequence, the solution can be obtained with a truncated SVD $W \Sigma V^\ast$ of $U'$, with $W, V$ made of the $M$ left and right singular vectors with the largest singular values and $\Sigma$ is the associated $M\times M$ diagonal singular values matrix. We set $\Phi = \text{Re}(W)$, $\Psi = \text{Im}(W)$ the real, resp. imaginary, part of $W$.

To find an appropriate value of the intermediate reduced dimension $M$ in the PSD-AE-HNN method, we set an error threshold and apply a search algorithm, such as dichotomy, to find the minimum value $M$ such that the reconstruction error is less than this threshold.


\clearpage

\printbibliography

\end{document}